\mag=1200
\input amstex
\documentstyle{amsppt}
\NoBlackBoxes
\TagsOnRight
\tolerance=1000

\define\hathatA{\,\,\hat{\!\!\hat{A}}}

\define\Otimes{\operatornamewithlimits{\otimes}}
\define\Oplus{\operatornamewithlimits{\oplus}}

\font\largeit=cmmi10 scaled \magstep2
\define\lM{\text{\largeit m}}

\let\varGamma=\Gamma
\let\varLambda=\Lambda
\let\ge=\geqslant

\let\myl=\overline
\let\wht=\widehat

\let\ot=\otimes
\let\myw=\widetilde
\let\mov=\overset
\let\mum=\setminus
\let\bt=\boxtimes
\let\myt=\times

\define\db#1{D^b({#1})}
\define\sO{{\Cal O}}
\define\sD{{\Cal D}}
\define\sF{{\Cal F}}
\define\sE{{\Cal E}}
\define\sG{{\Cal G}}
\define\sL{{\Cal L}}
\define\sH{{\Cal H}}
\define\sP{{\Cal P}}
\define\sJ{{\Cal J}}
\define\sM{{\Cal M}}
\define\QQ{\Bbb Q}
\define\RR{\Bbb R}
\define\CC{\Bbb C}
\define\ZZ{\Bbb Z}
\define\fA{\frak A}
\define\fD{\frak D}

\define\bR{\bold R}
\define\bL{\bold L}
\define\End{\Cal{E}nd}
\define\EEnd{\operatorname{End}}
\define\coh{\operatorname{coh}}
\define\har{\operatorname{char}}
\define\id{\operatorname{id}}
\define\Id{\operatorname{Id}}
\define\HHom{\Cal{H}\!\italic{om}}
\define\Hom{\operatorname{Hom}}
\define\codim{\operatorname{codim}}
\define\Auteq{\operatorname{Auteq}}
\define\Pic{\operatorname{Pic}}
\define\Aut{\operatorname{Aut}}
\define\Ad{\operatorname{Ad}}
\define\SL{\operatorname{SL}}
\define\Spec{\operatorname{Spec}}
\define\Ker{\operatorname{Ker}}
\define\Coker{\operatorname{Coker}}
\define\Iso{\operatorname{Iso}}
\define\Autoeq{\operatorname{Autoeq}}
\define\Supp{\operatorname{Supp}}
\define\NS{\operatorname{NS}}
\define\Ext{\operatorname{Ext}}
\define\ord{\operatorname{ord}}
\redefine\Sp{\operatorname{Sp}}
\define\sign{\operatorname{sign}}
\define\btl{\operatornamewithlimits{\bt}}
\redefine\Im{\operatorname{Im}}
\define\<{\langle}
\define\>{\rangle}
\topmatter
\date 1.X.2001\enddate
\translator A. V. Domrin \endtranslator
\address
V. A. Steklov Mathematical Institute of RAS
\endaddress
\email
orlov\@mi.ras.ru
\endemail
\author
D. O. Orlov
\endauthor
\title
Derived categories of coherent sheaves on abelian varieties and
equivalences between them
\endtitle
\abstract
We study derived categories of coherent sheaves on abelian varieties.
We give a criterion for the equivalence of the derived categories on
two abelian varieties. We describe the autoequivalence group for the
derived category of coherent sheaves of an abelian variety.
\endabstract
\endtopmatter
\rightheadtext{Derived categories of coherent sheaves}
\footnote""{This work is done under partial financial support of
the Russian Foundation for Basic Research (grant no.~99-01-01144),
grant in support of leading scientific schools
(no.~00-\allowbreak15-96085), and American Foundation of Civil Research
and Development for countries of the former Soviet Union
(grant no.~RM1-2089).}

\document

\head
Introduction
\endhead

For every algebraic variety~$X$ we have the abelian category~$\coh(X)$
of coherent sheaves on~$X$. Morphisms between varieties induce the
inverse image functors between the abelian categories of coherent sheaves.
If the morphism is proper, then the direct image functor is also defined.
These functors are not exact. They have left and right derived functors
respectively. To take all the derived functors into account, one must
pass from the abelian categories to their derived categories. For example,
to every smooth projective varieties there corresponds the so-called
bounded derived category of coherent sheaves $\db{X}$, and to every
morphism between such varieties there correspond the derived direct and
inverse image functors between the  derived categories of coherent sheaves.

The question is how many information is lost in the passage from varieties
to the derived categories of coherent sheaves. This passage actually
preserves almost all information. For example, in some cases one can
restore a variety from its derived category (see~\cite{2} or
Theorem~1.2 of this paper). Nevertheless, some classes of varieties
contain examples of two different varieties with equivalent derived
categories of coherent sheaves.

In this paper we study the case of abelian varieties.

Let $A$ be an abelian variety, and let $\hat A$ be the dual abelian variety.
As shown in~\cite{9}, the derived categories of coherent sheaves
$\db A$ and~$\db{\hat A}$ are equivalent and their equivalence, which is
called the {\it Fourier--Mukai transform}, can be defined in terms of the
Poincar\'e line bundle~$P_A$ on the product $A\myt\hat A$
by the following rule:
$$
F\mapsto\bR^\cdot p_{2*}\bigl(P_A\ot p_1^*(F)\bigr).
$$
This construction of Mukai was generalized in \cite{13} as follows.

Consider abelian varieties $A$,~$B$ and an isomorphism~$f$ between the
abelian varieties $A\myt\hat A$ and $B\myt\wht B$. We write~$f$ as a matrix
$$
\pmatrix
x & y
\\
z & w
\endpmatrix,
$$
where $x$ is a homomorphism from $A$ to $B$, \ $y$ is a homomorphism from
$\hat A$ to $B$ and so on. The isomorphism $f$ is called {\it isometric}
if its inverse map is given by
$$
f^{-1}=\pmatrix
\wht w & -\hat y
\\
-\hat z & \hat x
\endpmatrix.
$$

It is proved in \cite{13} that if $A$, $B$ are abelian varieties over an
algebraically closed field and there is an isometric isomorphism between
$A\myt\hat A$ and $B\myt\wht B$, then the derived categories of
coherent sheaves~$\db A$ and~$\db B$ are equivalent.

In this paper we prove the equivalence of these conditions over an
algebraically closed field of characteristic~$0$. In other words, we prove
that the derived categories $\db A$ and~$\db B$ are equivalent if and only
if there is an isometric isomorphism of $A\myt\hat A$ onto $B\myt\wht B$.
Part ``only if" of this assertion actually holds over an arbitrary field
(Theorem~2.19). As a corollary of this theorem, we obtain that for every
given abelian variety~$A$ there are only finitely many non-isomorphic
abelian varieties whose derived categories are isomorphic to $\db A$
(Corollary~2.20). The proof essentially uses the main theorem of~\cite{12},
which states that every exact equivalence between the derived categories
of coherent sheaves on smooth projective varieties can be represented by
an object on the product.

Representing equivalences by objects on the product, we construct a map
from the set of all exact equivalences between $\db A$ and~$\db B$
to the set of isometric isomorphisms of $A\myt\hat A$ onto $B\myt\wht B$.
Then we show that this map is functorial (Proposition~2.15). In particular,
we get a homomorphism from the group of exact autoequivalences of
$\db A$ to the group $U(A\myt\hat A)$ of isometric automorphisms of
$A\myt\hat A$.

The kernel of this homomorphism is described in \S\,3. It is isomorphic to
the direct sum of the free abelian group~$\ZZ$ and the group of $k$-points
of the variety $A\myt\hat A$ (Proposition~3.3). This fact is technically
based on Proposition~3.2, which states that the object on the product of
abelian varieties which determines an equivalence is actually a sheaf up
to a shift in the derived category. We note that this result holds for
abelian varieties only (for example, it is not true for K3-surfaces)
and is a key to the description of the autoequivalence group in the
case of abelian varieties.

In \S\,4 we assume that the ground field is algebraically closed and
$\har(k)=0$ and give another proof of the results of~\cite{13}. This proof
is based on the results of~\cite{10} that describe semihomogeneous bundles
on abelian varieties. In particular, we get an exact sequence of groups
$$
0\to\ZZ\oplus(A\myt\hat A)_k\to\Auteq\db A\to U(A\myt\hat A)\to1.
$$

In conclusion we describe the central extension of $U(A\myt\hat A)$ by
the group~$\ZZ$ and give a formula for the 2-cocycle that determines this
extension.

The author is grateful to the Max Planck Institute of Mathematics for
the hospitality and stimulating atmosphere. It is also a pleasure for
the author to express his gratitude to the Foundation for Promotion of
National Science.

\head
\S\,1. Preliminary information
\endhead

Let $X$ be an algebraic variety over a field~$k$ with structure sheaf
$\sO_X$. For every variety we have the abelian category~$\coh(X)$ of
coherent sheaves on it.

We denote by $\db X$ the bounded derived category of the abelian category
$\coh(X)$. It is obtained from the category of bounded complexes of
coherent sheaves by inversion of all quasi-isomorphisms, that is,
those maps between complexes that induce isomorphisms in cohomology
(see, for example,~\cite{3}).

Every derived category has the structure of a triangulated category.
This means that the additive category~$\sD$ is endowed with

a) an additive shift functor $[1]\:\sD\to\sD$, which is an
autoequivalence,

b) a class of distinguished triangles
$X\mov u\to\to Y\mov v\to\to Z\mov w\to\to X[1]$,
which must satisfy certain axioms (see~\cite{15}).

An additive functor $F\:\sD\to\sD'$ between triangulated categories is called
{\it exact} if it commutes with the shift functors and transforms every
distinguished triangle of~$\sD$ to a distinguished triangle of~$\sD'$.

In what follows we assume that all varieties are smooth and projective.
Any morphism $f\:X\to Y$ between smooth projective varieties induces two
exact functors: the direct image functor
$\bR^\cdot f_*\:\db X\to\db Y$ and the inverse image functor
$\bL^\cdot f^*\:\db Y\to\db X$. Moreover, each object $\sF\in\db X$
determines the exact functor of tensor product
$$
\Otimes^{\bL}\sF\:\db X\to\db X.
$$
Using these functors, we can introduce a large class of exact functors
between the categories~$\db X$ and~$\db Y$.

Let $X$, $Y$ be smooth projective varieties over the field~$k$.
We consider the Cartesian product~$X\myt Y$ and denote by
$p$,~$q$ the projections of $X\myt Y$ onto~$X$,~$Y$ respectively:
$$
X\mov p\to\leftarrow X\myt Y\mov q\to\to Y.
$$

Each object $\sE\in\db{X\myt Y}$ determines an exact functor $\Phi_{\sE}$
from the derived category~$\db X$ to the derived category $\db Y$ by the
formula
$$
\Phi_{\sE}(\cdot):=\bR^\cdot q_*\bigl(\sE\Otimes^{\bL}p^*(\cdot)\bigr).
\tag1.1
$$
Moreover, the same object $\sE\in\db{X\myt Y}$ determines another functor
$\Psi_{\sE}$ from the derived category~$\db Y$ to the derived category
$\db X$ by the formula similar to~\thetag{1.1}:
$$
\Psi_{\sE}(\cdot):=\bR^\cdot p_*\bigl(\sE\Otimes^{\bL}q^*(\cdot)\bigr).
$$
It is easily verified that $\Phi_{\sE}$ has left and right adjoint
functors~$\Phi_{\sE}^*$ and~$\Phi_{\sE}^!$ respectively. They
are given by
$$
\Phi_{\sE}^*\cong\Psi_{\sE^\vee\otimes q^*\omega_Y[\dim Y]}, \qquad
\Phi_{\sE}^!\cong\Psi_{\sE^\vee\ot p^*\omega_X[\dim X]}.
\tag1.2
$$
Here $\omega_X$ and $\omega_Y$ are the canonical sheaves on~$X$ and~$Y$,
and~$\sE^\vee$ is a convenient notation for
$\bR^\cdot\HHom(\sE,\sO_{X\myt Y})$.

Let now $X$, $Y$, $Z$ be three smooth projective varieties, and let
$\sE$, $\sF$, $\sG$ be objects of the derived categories $\db{X\myt Y}$,
$\db{Y\myt Z}$, $\db{X\myt Z}$ respectively. We consider the following
diagram of projections:


The objects $\sE$, $\sF$, $\sG$ determine the functors
$$
\align
&\Phi_{\sE}\:\db X\to\db Y,
\\
&\Phi_{\sF}\:\db Y\to\db Z,
\\
&\Phi_{\sG}\:\db X\to\db Z
\endalign
$$
by formula~\thetag{1.1}, that is,
$$
\align
\Phi_{\sE}
&:=\bR^\cdot\pi_{12*}^2\bigl(\sE\Otimes^{\bL}\pi_{12}^{1}{}^*(\cdot)\bigr),
\\
\Phi_{\sF}
&:=\bR^\cdot\pi^3_{23*}\bigl(\sF\Otimes^{\bL}\pi_{23}^{2}{}^*(\cdot)\bigr),
\\
\Phi_{\sG}
&:=\bR^\cdot\pi^3_{13*}\bigl(\sG\Otimes^{\bL}\pi_{13}^{1}{}^*(\cdot)\bigr).
\endalign
$$

We consider the object $p_{12}^*\sE\otimes p_{23}^*\sF\in\db{X\myt Y\myt Z}$.
In what follows we always denote it by $\sE\btl_Y\sF$. The following
assertion from \cite{9} yields the composition law for those exact functors
between derived categories that are represented by objects on the product.

\proclaim{Proposition 1.1}
The composite functor $\Phi_{\sF}\circ\Phi_{\sE}$ is isomorphic to
the functor~$\Phi_{\sG}$ represented by the object
$$
\sG=\bR^\cdot p_{13*}\bigl(\sE\btl_Y\sF\bigr).
\tag1.3
$$
\endproclaim

Thus, to every smooth projective variety there corresponds the derived
category of coherent sheaves on it, and to every object
$\sE\in\db{X\myt Y}$ on the product of such varieties there corresponds
the exact functor~$\Phi_{\sE}$ from the triangulated category $\db X$ to
the triangulated category~$\db Y$ with the composition law
described above.

The following two questions are fundamental in understanding this
correspondence.

1) When are derived categories of coherent sheaves on two different
smooth projective varieties equivalent as triangulated categories?

2) What is the group of exact autoequivalences of the
derived category of coherent sheaves on a given variety~$X$?

Some results in this direction are already known. For example, there is
a complete answer to these questions in the case when either canonical or
anticanonical sheaf of the variety is ample.

\proclaim{Theorem 1.2~\cite{2}}
Let $X$ be a smooth projective variety whose canonical (or anticanonical)
sheaf is ample. Suppose that the category $\db X$ is equivalent (as a
triangulated category) to the derived category~$\db{X'}$ for some smooth
algebraic variety~$X'$. Then $X'$ is isomorphic to~$X$.
\endproclaim

\proclaim{Theorem 1.3~\cite{2}}
Let $X$ be a smooth projective variety whose canonical (or anticanonical)
sheaf is ample. Then the group of isomorphism classes of exact
autoequivalences of the category $\db X$ is generated by automorphisms
of the variety, twists by line bundles and shifts in the derived
category.
\endproclaim

The group of exact autoequivalences may also be described. For any
variety~$X$, the group $\Auteq\db X$ of exact autoequivalences always
contains a subgroup $G(X)$ which is a semidirect product of its normal
subgroup $G_1=\Pic(X)\oplus\ZZ$ and the subgroup $G_2=\Aut X$ with its
natural action on~$G_1$. Under the inclusion $G(X)\subset\Auteq\db X$,
the generator of~$\ZZ$ is mapped to the shift functor~$[1]$,
each line bundle $\sL\in\Pic(X)$ is mapped to the functor~$\ot\sL$,
and each automorphism
$f\:X\to X$ induces the autoequivalence~$\bR^\cdot f_*$.

Under the hypotheses of Theorem~1.3, we can additionally assert that
the group $\Auteq\db X$ of exact autoequivalences coincides with
$G(X)$, that is,
$$
\Auteq\db X\cong\Aut X\ltimes\bigl(\Pic(X)\oplus\ZZ\bigr).
$$

To study the cases of equivalence of the derived categories of coherent
sheaves on two varieties and to describe their autoequivalence groups,
it is desirable to have explicit formulae for all exact functors.
It is conjectured that they are all represented by objects on the
product, that is, are given by~\thetag{1.1}. This conjecture is
proved in the particular case of equivalences.

\proclaim{Theorem 1.4~\cite{12}}
Let $X$, $Y$ be smooth projective varieties. Suppose that
$F\:\db X\mov\sim\to\to\db Y$ is an exact functor and an equivalence
of triangulated categories. Then there is a unique (up to an isomorphism)
object $\sE\in\db{X\myt Y}$ such that the functor $F$ is isomorphic to
the functor~$\Phi_{\sE}$.
\endproclaim

To verify that a functor $F$ is an equivalence, it suffices to show that
$F$ and its right (or left) adjoint functor are fully faithful. We
recall that a functor~$F$ is  {\it fully faithful} if, for any objects
$A$ and~$B$, the natural map
$$
\Hom(A,B)\to\Hom\bigl(F(A),F(B)\bigr)
$$
is a bijection. In what follows we need a method to determine whether
the functors $\Phi_{\sE}\:\db X\to\db Y$ are fully faithful. This is
rather difficult to verify in general, but the following criterion is
useful in some situations.

\proclaim{Theorem 1.5~\cite{1}}
Let $M$, $X$ be smooth projective varieties over an algebraically closed
field of characteristic~$0$. Take $\sE\in\db{M\myt X}$. Then the functor
$\Phi_{\sE}$ is fully faithful if and only if the following orthogonality
conditions hold:

\rom{1)} $\Hom^i_X\bigl(\Phi_{\sE}(\sO_{t_1}),\Phi_{\sE}(\sO_{t_2})\bigr)=0$
for all~$i$ and $t_1\ne t_2$;

\rom{2)} $\Hom^0_X\bigl(\Phi_{\sE}(\sO_t),\Phi_{\sE}(\sO_t)\bigr)=k$;

\rom{3)} $\Hom^i_X\bigl(\Phi_{\sE}(\sO_t),\Phi_{\sE}(\sO_t)\bigr)=0$ for
$i\notin[0,\dim M]$.

Here $t$, $t_1$, $t_2$ are points of $M$, and $\sO_{t_i}$ are the
corresponding skyscraper sheaves.
\endproclaim

Let now $X_1$, $X_2$, $Y_1$, $Y_2$ be smooth projective varieties.
We consider objects~$\sE_1$ and~$\sE_2$ of
the categories $\db{X_1\myt Y_1}$ and $\db{X_2\myt Y_2}$ respectively.
By definition, the object
$$
\sE_1\bt\sE_2\in\db{(X_1\myt X_2)\myt(Y_1\myt Y_2)},
$$
is equal to $p^*_{13}(\sE_1)\otimes p^*_{24}(\sE_2)$. As before
(see~\thetag{1.1}), the objects $\sE_1$, $\sE_2$, $\sE_1\bt\sE_2$
determine the functors
$$
\align
\Phi_{\sE_1}\:\db{X_1}
&\to\db{Y_1},
\\
\Phi_{\sE_2}\:\db{X_2}
&\to\db{Y_2},
\\
\Phi_{\sE_1\bt\sE_2}\:\db{X_1\myt X_2}
&\to\db{Y_1\myt Y_2}.
\endalign
$$

We consider an object $\sG\in\db{X_1\myt X_2}$ and denote the object
$\Phi_{\sE_1\bt\sE_2}(\sG)\in\db{Y_1\myt Y_2}$ by~$\sH$. By~\thetag{1.1},
these objects determine the functors
$$
\Phi_{\sG}\:\db{X_1}\to\db{X_2}, \qquad
\Phi_{\sH}\:\db{Y_1}\to\db{Y_2}.
$$

\proclaim{Lemma 1.6}
In the notation above, we have an isomorphism of functors
$\Phi_{\sH}\cong\Phi_{\sE_2}\circ\Phi_{\sG}\circ\Psi_{\sE_1}$.
\endproclaim

\demo\nofrills {The proof} \  immediately follows from Proposition~1.1.
\enddemo

Let $Z_1$, $Z_2$ be further smooth projective varieties, and let
$\sF_1$, $\sF_2$ be objects of the categories $\db{Y_1\myt Z_1}$,
$\db{Y_2\myt Z_2}$ respectively.
Then we also have the functors $\Phi_{\sF_1}$,
$\Phi_{\sF_2}$, $\Phi_{\sF_1\bt\sF_2}$.
By~\thetag{1.3}, we can find objects $\sG_1$ and~$\sG_2$
of the categories $\db{X_1\myt Z_1}$ and $\db{X_2\myt Z_2}$ such that
$$
\Phi_{\sG_1}\cong\Phi_{\sF_1}\circ\Phi_{\sE_1}, \qquad
\Phi_{\sG_2}\cong\Phi_{\sF_2}\circ\Phi_{\sE_2}.
$$
It is directly verified that there is a natural relation
$$
\Phi_{\sF_1\bt\sF_2}\circ\Phi_{\sE_1\bt\sE_2}\cong\Phi_{\sG_1\bt\sG_2}.
\tag1.4
$$
Using this relation, we easily prove the following assertion.

\proclaim{Assertion 1.7}
Under the conditions above, assume that the functors
$\Phi_{\sE_1}$ and~$\Phi_{\sE_2}$ are fully faithful
(resp. equivalences). Then the functor
$$
\Phi_{\sE_1\bt\sE_2}\:\db{X_1\myt X_2}\to\db{Y_1\myt Y_2}
$$
is also fully faithful (resp. an equivalence).
\endproclaim

\demo{Proof}
Let $F$ be a functor having an adjoint (say, the left adjoint
$F^*$). Then $F$ is fully faithful if and only if the composite
$F^*\circ F$ is isomorphic to the identity functor. The functors
$\Phi_{\sE_i}$ have left adjoints $\Phi_{\sE_i}^*$, which are
defined by~\thetag{1.2}. If they are fully faithful, then
the composites $\Phi_{\sE_i}^*\circ\Phi_{\sE_i}$ are isomorphic to
the identity functors, which are presented by the structure sheaves
of the diagonals $\Delta_i\subset X_i\myt X_i$. We easily verify that
the sheaf $\sO_{\Delta_1}\bt\sO_{\Delta_2}$ is isomorphic to the structure
sheaf~$\sO_\Delta$, where $\Delta$ is the diagonal in
$(X_1\myt X_2)\myt(X_1\myt X_2)$.

Using~\thetag{1.4}, we see that the composite
$\Phi_{\sE_1\bt\sE_2}^*\circ\Phi_{\sE_1\bt\sE_2}$ is represented by the
structure sheaf of the diagonal~$\Delta$. Hence it is isomorphic to the
identity functor. Thus $\Phi_{\sE_1\bt\sE_2}$ is fully faithful.
The assertion on the equivalences is proved similarly.
\enddemo

Suppose that the functor $\Phi_{\sE}\:\db X\to\db Y$ is an equivalence and
that an object $\sF\in\db{X\myt Y}$ satisfies
$\Psi_{\sF}\cong\Phi_{\sE}^{-1}$. Then we denote the functor
$$
\Phi_{\sF\bt\sE}\:\db{X\myt X}\to\db{Y\myt Y}
\tag1.5
$$
by $\Ad_{\sE}$. The functor $\Ad_{\sE}$ is an equivalence by Assertion~1.7.
Moreover, by Lemma~1.6, for every object $\sG\in\db{X\myt X}$ there is
an isomorphism of functors
$$
\Phi_{\Ad_{\sE}(\sG)}\cong\Phi_{\sE}\circ\Phi_{\sG}\circ\Phi_{\sE}^{-1}.
\tag1.6
$$

\head
\S\,2. Equivalences between the categories of coherent sheaves on abelian
varieties
\endhead

Let $A$ be an abelian variety of dimension $n$ over a field~$k$. We denote
by $m\:A\myt A\to A$ the composition morphism (regarded as a morphism
defined over~$k$) and let $e$ be the $k$-point which is the identity of the
group structure.

We denote by $\hat A$ the dual abelian variety, which is the moduli space
of line bundles on~$A$ that belong to~$\Pic^0(A)$. Moreover, $\hat A$ is
a fine moduli space. Hence the product $A\myt\hat A$ carries the
universal line bundle~$P$, which is called the {\it Poincar\'e bundle}.
This bundle is uniquely determined by the condition that, for every
$k$-point $\alpha\in\hat A$, the restriction of~$P$ to $A\myt\{\alpha\}$
is isomorphic to the line bundle in~$\Pic^0(A)$ that corresponds
to~$\alpha$, and the restriction $P|_{\{e\}\myt\hat A}$ is trivial.

\definition{Definition 2.1}
In what follows we denote by $P_\alpha$ the line bundle on~$A$
that corresponds to the $k$-point $\alpha\in\hat A$.
Moreover, suppose that $A_1,\dots,A_m$ are abelian varieties,
and let $(\alpha_1,\dots,\alpha_m)$ be a $k$-point of
$\hat A_1\myt\cdots\myt\hat A_m$. Then we denote by
$P_{(\alpha_1,\dots,\alpha_k)}$ the line bundle
$P_{\alpha_1}\bt\cdots\bt P_{\alpha_k}$ on the product
$A_1\myt\cdots\myt A_k$.
\enddefinition

For any homomorphism $f\:A\to B$ of abelian varieties we have the dual
homomorphism $\hat f\:\wht B\to\hat A$. It maps a point $\beta\in\wht B$
to a point $\alpha\in\hat A$ if and only if the line bundle
$f^*P_\beta$ coincides with the bundle~$P_\alpha$ on~$A$.

The second dual abelian variety~$\hathatA$ can naturally be
identified with~$A$ using the Poincar\'e bundles on $A\myt\hat A$ and
$\hat A\myt\hathatA$. In other words, there is a unique isomorphism
$\kappa_A\:A\mov\sim\to\to\hathatA$ such that the lifting of the Poincar\'e
bundle~$P_{\hat A}$ by the isomorphism
$1\myt\kappa_A\:\hat A\myt A\mov\sim\to\to\hat A\myt\hathatA$
coincides with the Poincar\'e bundle~$P_A$, that is,
$(1\myt\kappa_A)^*P_{\hat A}\cong P_A$.

Hence the correspondence $\wht{\,}$ is an involution on the category of
abelian varieties, that is, $\wht{\,}$ is a contravariant functor whose
square is isomorphic to the identity functor. The isomorphism is given
by the map $\kappa\:\Id\mov\sim\to\to\wht{\wht{\,}}$\,\,.

\remark{Remark \rm2.2 ($k=\CC$)}
Let $A$ be an abelian variety over~$\CC$. We choose a basis
$l_1,\dots,l_{2n}$ in~$H_1(A,\ZZ)$, and let $l_1^*,\dots,l^*_{2n}$ be the
dual basis in $H_1(\hat A,\ZZ)\cong H_1(A,\ZZ)^*$. Let
$l^{**}_1,\dots,l^{**}_{2n}$ be the basis in $H_1(\hathatA,\ZZ)$
dual to $l^*_1,\dots,l^*_{2n}$.
The isomorphism $\kappa_A\:A\mov\sim\to\to\hathatA$ induces an identification
of~$H_1(A,\ZZ)$ with~$H_1(\hathatA,\ZZ)$ such that the elements~$l_i$ are
identified with~$-l^{**}_i$. The sign ``minus'' appears because the forms
$c_1(P_A)$ and~$c_1(P_{\hat A})$ are skew-symmetric.
\endremark

\remark{Remark \rm2.3 ($k=\CC$)}
Let $f\:A\to B$ be a homomorphism of complex abelian varieties, and let
$\hat f\:\hat A\to\wht B$ be the dual homomorphism. We fix bases in the
homology groups $H_1(A,\ZZ)$,~$H_1(B,\ZZ)$ and take the dual bases in the
first homology groups of $\hat A$,~$\wht B$. The maps $f$ and~$\hat f$
induce linear maps between the first homology groups, and we denote
their matrices by $F$ and~$\wht F$. Then the matrix~$\wht F$ is transposed
to~$F$.

Consider a homomorphism $f\:A\to\hat A$. By the definition of the
isomorphism $\kappa$, we may assume that~$\hat f$ is also a homomorphism
from~$A$ to~$\hat A$. Hence the matrices $F$ and~$\wht F$ are
skew-transposed to each other, that is, $\wht F=-F^t$.
\endremark

The Poincar\'e bundle $P$ provides an example of an exact equivalence
between the derived categories of coherent sheaves on two (generally
non-isomorphic) varieties $A$ and~$\hat A$. We consider the projections
$$
A\mov p\to\leftarrow A\myt\hat A\mov q\to\to\hat A
$$
and define a functor $\Phi_P\:\db A\to\db{\hat A}$ by~\thetag{1.1},
that is, $\Phi_P(\cdot)=\bR q_*\bigl(P\ot p^*(\cdot)\bigr)$.

\proclaim{Proposition 2.4~\cite{9}}
Let $P$ be the Poincar\'e bundle on $A\myt\hat A$. Then the functor
$\Phi_P\:\db A\to\db{\hat A}$ is an exact equivalence, and there is an
isomorphism of functors
$$
\Psi_P\circ\Phi_P\cong(-1_A)^*[n],
$$
where $(-1_A)$ is the morphism of group inversion.
\endproclaim

\remark{Remark \rm 2.5}
This assertion is proved in~\cite{9} for abelian varieties over an
algebraically closed field. However, it holds over an arbitrary field
since the dual variety and the Poincar\'e bundle are always defined
over the same field (see~\cite{11}, for example) while the assertion
on the equivalence of the functor will follow from Lemma~2.12 (see~below).
\endremark

For every $k$-point $a\in A$ there is the shift automorphism
$m(\cdot,a)\:A\to A$, which will be denoted by~$T_a$. For every $k$-point
$\alpha\in\hat A$, we denote by~$P_\alpha$ the corresponding line bundle
on~$A$.

We now consider a $k$-point $(a,\alpha)\in A\myt\hat A$. It determines a
functor from~$\db A$ to itself by
$$
\Phi_{(a,\alpha)}(\cdot):=T_{a*}(\cdot)\ot P_\alpha
=T^*_{-a}(\cdot)\ot P_\alpha.
\tag2.1
$$
The functor $\Phi_{(a,\alpha)}$ is represented by the sheaf
$$
S_{(a,\alpha)}=\sO_{\Gamma_a}\ot p_2^*(P_\alpha)
\tag2.2
$$
on the product $A\myt A$, where $\Gamma_a$ is the graph of the shift
automorphism~$T_a$. The functor~$\Phi_{(a,\alpha)}$ is clearly an
autoequivalence.

The set of functors $\Phi_{(a,\alpha)}$ (parametrized by
$A\myt\hat A$) can be combined into a single functor from $\db{A\myt\hat A}$ to
$\db{A\myt A}$, which maps the structure sheaf $\sO_{(a,\alpha)}$
of the point to~$S_{(a,\alpha)}$.
(We note that this condition determines the functor
non-uniquely but only up to multiplying by a line bundle lifted from
$A\myt\hat A$.)

We define a functor $\Phi_{S_A}\:\db{A\myt\hat A}\to\db{A\myt A}$
as a composite of two other functors. To do this, we consider the object
$\sP_A=p^*_{14}\sO_\Delta\ot p^*_{23}P\in\db{(A\myt\hat A)\myt(A\myt A)}$
and denote by $\mu_A\:A\myt A\to A\myt A$ the morphism that sends a point
$(a_1,a_2)$ to $\bigl(a_1,m(a_1, a_2)\bigr)$. We get two functors:
$$
\Phi_{\sP_A}\:D^b(A\myt\hat A)\to D^b(A\myt A), \qquad
\bR\mu_{A*}\:\db{A\myt A}\to\db{A\myt A}.
$$

\definition{Definition 2.6}
The functor $\Phi_{S_A}$ is the composite $\bR\mu_{A*}\circ\Phi_{\sP_A}$.
\enddefinition

We can explicitly describe the object $S_A$ on the product
$(A\myt\hat A)\myt(A\myt A)$ that represents the functor~$\Phi_{S_A}$.
Since the explicit formula is not used in what follows,
we present it without proof.

\proclaim{Lemma 2.7}
Let $S_A$ be the object on the product $(A\myt\hat A)\myt(A\myt A)$
that represents the functor~$\Phi_{S_A}$. Then
$$
S_A=(m\cdot p_{13},p_4)^*\sO_\Delta\ot p_{23}^*P_A.
$$
Here $(m\cdot p_{13},p_4)$ is the morphism onto $A\myt A$ that maps
a point $(a_1,\alpha,a_3,a_4)$ to the point
$\bigl(m(a_1,a_3),a_4\bigr)$.
\endproclaim

\proclaim{Assertion 2.8}
The functor $\Phi_{S_A}$ is an equivalence and, for every $k$-point
$(a,\alpha)\in A\myt\hat A$,

{\rm a)} $\Phi_{S_A}$ maps the structure sheaf $\sO_{(a,\alpha)}$ of the
point to the sheaf~$S_{(a,\alpha)}$ defined by~\thetag{2.2},

{\rm b)} $\Phi_{S_A}$ maps the line bundle $P_{(\alpha,a)}$ on
$A\myt\hat A$ to the object $\sO_{\{-a\}\myt A}\otimes p_2^* P_\alpha[n]$.
\endproclaim

\demo{Proof}
By definition, $\Phi_{S_A}$ is the composite of the functors $\bR\mu_{A*}$
and~$\Phi_{\sP_A}$, which are equivalences. (This is obvious for the first
functor, and this follows from Assertion~1.7 and Proposition~2.4 for the
second one.)

The functor $\Phi_{\sP_A}$ maps the structure sheaf~$\sO_{(a,\alpha)}$
 to the sheaf~$\sO_{A\myt\{a\}}\ot p_1^*P_\alpha$, and the
functor~$\bR\mu_{A*}$ maps the sheaf
$\sO_{A\myt\{a\}}\ot p_1^*P_\alpha$ to the sheaf
$\sO_{\Gamma_a}\ot p_1^*(P_\alpha)$.

Using Proposition~2.4, we find in the same way that the functor
$\Phi_{\sP_A}$ maps the line bundle~$P_{(\alpha,a)}$ to the object
$\sO_{\{-a\}\myt A}\ot p_2^*P_\alpha[n]$, and the functor~$\bR\mu_{A*}$
maps the object $\sO_{\{-a\}\myt A}\ot p_2^*P_\alpha[n]$ to itself.
\enddemo

Suppose that $A$ and $B$ are abelian varieties whose
derived categories of coherent sheaves are equivalent.
We fix an equivalence. By Theorem~1.4, it is represented by an object
on the product. Thus we have an object
$\sE\in\db{A\myt B}$ and the equivalence
$\Phi_{\sE}\:\db A\mov\sim\to\to\db B$.

We recall that the functor
$$
\Ad_{\sE}\:\db{A\myt A}\mov\sim\to\to\db{B\myt B}
$$
is defined by~\thetag{1.5} and is an equivalence. We consider
the composite functor $\Phi_{S_B}^{-1}\circ\Ad_{\sE}\circ\Phi_{S_A}$.

\definition{Definition 2.9}
We denote by $\sJ(\sE)$ the object that represents the functor
$$
\Phi_{S_B}^{-1}\circ\Ad_{\sE}\circ\Phi_{S_A}.
$$
\enddefinition

Thus there is a commutative diagram
$$
\CD
\db{A\myt\hat A}@>\Phi_{S_A}>>\db{A\myt A}
\\
@V \Phi_{\!\sJ(\sE)}VV @VV\Ad_{\sE}V
\\
\db{B\myt\wht B}@>\Phi_{S_B}>>\db{B\myt B}
\endCD
\tag2.3
$$

The following theorem enables us to describe the object $\sJ(\sE)$
and is basic for the description of abelian varieties with equivalent
derived categories of coherent sheaves.

\proclaim{Theorem 2.10}
One can find a homomorphism $f_{\sE}\:A\myt\hat A\to B\myt\wht B$
of abelian varieties and a line bundle~$L_{\sE}$ on $A\myt\hat A$ such that
$f_{\sE}$ is an isomorphism and the object~$\sJ(\sE)$ is isomorphic to
$i_*(L_{\sE})$, where~$i$ is the embedding of $A\myt\hat A$ into
$(A\myt\hat A)\myt(B\myt\wht B)$ as the graph of the isomorphism~$f_{\sE}$.
\endproclaim

Before proving this theorem, we state two lemmas that enable us to assume
that $k$ is algebraically closed. We denote the algebraic closure of~$k$ by
$\myl k$. We put $\myl X:=X\myt_{\Spec(k)}\Spec(\myl k)$ and denote by
$\myl\sF$ the inverse image of~$\sF$ under the morphism
$\myl X\to X$.

\proclaim{Lemma 2.11}
Let $\sF$ be a coherent sheaf on a smooth variety~$X$. Suppose that we have
a closed subvariety $j\:Z\hookrightarrow\myl X$ and an invertible sheaf
$\sL$ on~$Z$ such that $\myl\sF\cong j_*\sL$. Then we can find a closed
subvariety $i\:Y\hookrightarrow X$ and an invertible sheaf $\sM$ on~$Y$
such that $\sF\cong i_*\sM$ and $j=\myl i$.
\endproclaim

\demo{Proof} Since the argument is local, we may assume that
we have an affine variety $X=\Spec(A)$ and an $A$-module~$M$.
Let $J\subset A$ be the annihilator of the module~$M$, and let
$J'\subset\myl A=A\otimes_k\myl k$ be the annihilator of the module
$\myl M=M\otimes_k\myl k$. If $\{e_i\}$ is a basis of the field~$\myl k$
over~$k$, then $\myl M=\oplus Me_i$ as a module over~$A$. Clearly,
$J\otimes_k\myl k\subseteq J'$. On the other hand, if an element
$\sum a_i\otimes e_i\in\myl A$ belongs to~$J'$, then $\sum a_im\otimes e_i=0$
for every $m\in M$. Hence each~$a_i$ belongs to~$J$. Thus,
$J\otimes_k\myl k=J'$.

The hypothesis of the lemma implies that $\myl M$ is a projective module
of rank~1 over the algebra $\myl B:=\myl A/J'=B\otimes_k\myl k$, where
$B=A/J$, and $\myl M=M\otimes_B\myl B$. Since~$\myl B$ is a strictly flat
$B$-algebra, we see that~$M$ is a projective $B$-module of rank~1 (see,
for example, \cite{4},~\thetag{3.1.4}).
\enddemo

The following lemma asserts that the property of a functor to be fully
faithul (or an equivalence) is stable with respect to field extensions.

\proclaim{Lemma 2.12}
Let $X$, $Y$ be smooth projective varieties over~$k$, and let $\sE$
be an object of the derived category $\db{X\myt Y}$. Consider a field
extension $F/k$ and the varieties
$$
X'=X\underset{\Spec(k)}\to\myt\Spec(F), \qquad
Y'=Y\underset{\Spec(k)}\to\myt\Spec(F).
$$
Let $\sE'$ be the lifting of the object $\sE$ to the category
$\db{X'\myt Y'}$. The functor $\Phi_{\sE}\:\db X\to\db Y$ is fully
faithful (resp. an equivalence) if and only if the functor
$\Phi_{\sE'}\:\db{X'}\to\db{Y'}$ is fully faithful (resp. an
equivalence).
\endproclaim

\demo{Proof}
As above, we denote by~$\Phi_{\sE}^*$ the left adjoint to the functor
$\Phi_{\sE}$. If $\Phi_{\sE}$ is fully faithful, then the composite
$\Phi_{\sE}\circ\Phi_{\sE}^*$ is the identity functor $\id_{\db X}$,
which is known to be represented by the structure sheaf $\sO_\Delta$
of the diagonal in the product $X\myt X$. Using Proposition~1.1
and the theorem on flat base change, we see that the composite
$\Phi_{\sE'}\circ\Phi_{\sE'}^*$ is represented by the structure sheaf
$\sO_{\Delta'}$, where~$\Delta'$ is the diagonal in $X'\myt X'$. Hence
the functor $\Phi_{\sE'}$ is fully faithful.

Conversely, consider the composite $\Phi_{\sE}\circ\Phi_{\sE}^*$. It is
represented by some object~$\sJ$ on $X\myt X$. There is a canonical
morphism $\phi\:\sJ\to\sO_\Delta$. Since $\Phi_{\sE'}$ is fully faithful
by the hypothesis, the morphism $\phi'\:\sJ'\to\sO_{\Delta'}$ is an
isomorphism. It follows immediately that~$\phi$ is an isomorphism as well,
whence the functor~$\Phi_{\sE}$ is fully faithful.

We similarly prove the assertion on equivalences, which follows
since the adjoint functor is fully faithful.
\enddemo

\demo{Proof of Theorem 2.10}
Using Lemmas 2.11 and~2.12, we may pass to the algebraic closure of~$k$.

{\it Step\/} 1. We denote by $e$ and $e'$ the closed points of
$A\myt\hat A$ and $B\myt\wht B$ (respectively) that are the identity
elements of the group structures. We consider the skyscraper sheaf $\sO_e$
and calculate its image under the functor~$\Phi_{\!\sJ(\sE)}$. By
definition, we know that
$$
\Phi_{\!\sJ(\sE)}=\Phi_{S_B}^{-1}\circ\Ad_{\sE}\circ\Phi_{S_A}.
$$
By Assertion 2.8, the functor $\Phi_{S_A}$ maps the sheaf $\sO_e$ to the
structure sheaf $\sO_{\Delta(A)}$ of the diagonal in $A\myt A$. Since the
structure sheaf of the diagonal represents the identity functor,
\thetag{1.6} implies that $\Ad_{\sE}(\sO_{\Delta(A)})$ is the structure
sheaf $\sO_{\Delta(B)}$ of the diagonal in $B\myt B$. By Assertion~2.8,
the functor~$\Phi_{S_B}^{-1}$ maps the last sheaf to the structure
sheaf~$\sO_{e'}$.

{\it Step\/} 2. We thus obtain that
$$
\sJ(\sE)\Oplus^{\bL}\sO_{\{e\}\myt(B\myt\wht B)}\cong\sO_{\{e\}\myt\{e'\}}.
$$
It follows that there is an affine neighbourhood $U=\Spec(R)$ of~$e$ in the
Zariski topology such that the object $\sJ':=\sJ(\sE)|_{U\myt(B\myt\wht B)}$
is a coherent sheaf whose support intersects the fibre
$\{e\}\myt(B\myt\wht B)$ at the point $\{e\}\myt\{e'\}$.
We recall that the support of a coherent sheaf is a closed subset.

We now consider an affine neighbourhood $V=\Spec(S)$ of~$e'$ in
$B\myt\wht B$. The intersection of the support of $\sJ'$ with the
complement $B\myt\wht B\mum V$ is a closed subset whose projection
onto $A\myt\hat A$ is a closed subset that does not contain~$e$.

Shrinking $U$ if necessary, we may thus assume that it is still affine
and the support of~$\sJ'$ is contained in $U\myt V$. This means that
there is a coherent sheaf~$\sF$ on $U\myt V$ such that $j_*(\sF)=\sJ'$,
where~$j$ is the inclusion of $U\myt V$ into $U\myt(B\myt\wht B)$.
We denote by~$M$ the finitely generated $R\otimes S$-module corresponding
to the sheaf~$\sF$, that is, $\sF=\myw M$. Moreover, $M$ is a finitely
generated $R$-module since the direct image of the coherent sheaf
$\sJ'=j_*\sF$ under the projection is a coherent sheaf.

Let $m$ be the maximal ideal in $R$ corresponding to the point~$e$.
It is known that
$$
M\otimes_R R/m\cong R/m.
$$
Hence there is a homomorphism $\phi\:R\to M$ of $R$-modules, which becomes
an isomorphism after the tensor multiplication by~$R/m$. Hence the supports
of the coherent sheaves $\Ker\phi$ and~$\Coker\phi$ do not contain~$e$.
Replacing~$U$ by a smaller affine neighbourhood of~$e$, disjoint from
the supports of $\Ker\phi$ and~$\Coker\phi$, we see that~$\phi$ is an
isomorphism. Hence there is a subscheme $X(U)\subset U\myt(B\myt\wht B)$
such that the projection $X(U)\to U$ is an isomorphism and
$$
\sJ'=\sJ(\sE)|_{U\myt(B\myt\wht B)}\cong\sO_{X(U)}.
$$

{\it Step\/} 3. We thus see that for every closed point
$(a,\alpha)\in U$ there is a closed point
$(b,\beta)\in B\myt\wht B$ such that
$$
\Phi_{\!\sJ(\sE)}(\sO_{(a,\alpha)})\cong\sO_{(b,\beta)}.
$$
Any closed point $(a,\alpha)\in A\myt\hat A$ may be presented as a sum
$(a,a')=(a_1,\alpha_1)+(a_2,\alpha_2)$, where the points
$(a_1,\alpha_1)$, $(a_2,\alpha_2)$ belong to~$U$. We denote by
$(b_1,\beta_1)$ and $(b_2,\beta_2)$ the images of these points under the
functor $\Phi_{\!\sJ(\sE)}$. It is known that the functor $\Phi_{S_A}$
maps the structure sheaf $\sO_{(a,\alpha)}$ to the sheaf $S_{(a,\alpha)}$.
We denote the object $\Ad_{\sE}(S_{(a,\alpha)})$ by~$\sG$. To calculate it,
we use~\thetag{1.6}. We have
$$
\Phi_{\sG}\cong\Phi_{\sE}\circ\Phi_{(a,\alpha)}\circ\Phi_{\sE}^{-1}.
$$
But the functor $\Phi_{(a,\alpha)}$, which equals
$T_a^*(\cdot)\ot P_\alpha$ by definition~\thetag{2.1},
is represented as the composite
$\Phi_{(a_1,\alpha_1)}\circ\Phi_{(a_2,\alpha_2)}$.
We thus get a sequence of isomorphisms
$$
\alignat1
\Phi_{\sG}
&\cong\Phi_{\sE}\circ\Phi_{(a,\alpha)}\circ\Phi_{\sE}^{-1}
\cong\Phi_{\sE}\circ\Phi_{(a_1,\alpha_1)}\circ\Phi_{\sE}^{-1}
\cong\Phi_{\sE}\circ\Phi_{(a_2,\alpha_2)}\circ\Phi_{\sE}^{-1}
\\
&\cong\Phi_{(b_1,\beta_1)}\circ\Phi_{(b_2,\beta_2)}\cong\Phi_{(b,\beta)},
\endalignat
$$
where $(b,\beta)=(b_1,\beta_1)+(b_2,\beta_2)$. Hence the object $\sG$ is
isomorphic to $S_{(b,\beta)}$. We finally obtain that
$$
\Phi_{\!\sJ(\sE)}(\sO_{(a,\alpha)})\cong\sO_{(b,\beta)}
$$
for every closed point $(a,\alpha)\in A\myt\hat A$.

For every closed point $(a,a')$, we can now repeat
the procedure of Step 2 to find a neighbourhood~$W$ and a subscheme
$X(W)\subset W\myt(B\myt\wht B)$ such that the projection $X(W)\to W$
is an isomorphism and $\sJ|_{W\myt(B\myt\wht B)}\cong\sO_{X(W)}$.

Gluing all these neighbourhoods, we find a subvariety
$i\:X\hookrightarrow(A\myt\hat A)\myt(B\myt\wht B)$
such that the projection $X\to A\myt\hat A$ is an isomorphism
and the sheaf~$\sJ(\sE)$ is isomorphic to the sheaf~$i_*L$,
where $L$ is a line bundle on~$X$. The subvariety~$X$ determines
a homomorphism from $A\myt\hat A$ to $B\myt\wht B$ which induces
an equivalence of derived categories. Hence this homomorphism is
an isomorphism.
\enddemo

In particular, Theorem 2.10 immediately implies that if abelian varieties
$A$ and $B$ have equivalent derived categories of coherent sheaves, then
the varieties $A\myt\hat A$ and $B\myt\wht B$ are isomorphic.
Below we shall see that this isomorphism must satisfy an additional condition
(see Proposition~2.18).

\proclaim{Corollary 2.13}
The isomorphism $f_{\sE}$ maps a $k$-point $(a,\alpha)\in A\myt\hat A$
to a point $(b,\beta)\in B\myt\wht B$ if and only if the equivalences
$$
\Phi_{(a,\alpha)}\:\db A\mov\sim\to\to\db A, \qquad
\Phi_{(b,\beta)}\:\db B\mov\sim\to\to\db B
$$
defined by~\thetag{2.1}, satisfy the relation
$$
\Phi_{(b,\beta)}\circ\Phi_{\sE}\cong\Phi_{\sE}\circ\Phi_{(a,\alpha)}.
$$
This relation is equivalent to the following condition on the objects that
represent these functors:
$$
T_{b*}\sE\ot P_\beta\cong T_{-a*}\sE\ot P_\alpha=T^*_a\sE\ot P_\alpha.
$$
\endproclaim

\demo{Proof}
By Theorem 2.10, the functor $\Phi_{\!\sJ(\sE)}$ maps the structure sheaf
$\sO_{(a,\alpha)}$ of the point $(a,\alpha)$
to the structure sheaf~$\sO_{(b,\beta)}$
of the point $(b,\beta)=f_{\sE}(a,\alpha)$.
Assertion~2.8 implies that~$\Phi_{S_A}$ maps the sheaf~$\sO_{(a,\alpha)}$
to~$S_{(a,\alpha)}$. The sheaf~$S_{(a,\alpha)}$ in its turn represents
the functor
$$
\Phi_{(a,\alpha)}=T_{a*}(\cdot)\ot P_\alpha.
$$
Using diagram~\thetag{2.3}, we now see that~$f_{\sE}$ maps the point
$(a,\alpha)$ to $(b,\beta)$ if and only if
$S_{(b,\beta)}\cong\Ad_{\sE}(S_{(a,\alpha)})$.
Applying~\thetag{1.6}, we find that
$\Phi_{(b,\beta)}\cong\Phi_{\sE}\circ\Phi_{(a,\alpha)}\circ\Phi_{\sE}^{-1}$.
\enddemo

In what follows we use an explicit formula for the object $\sJ(\sE)$
in the particular case when $A=B$ and the equivalence~$\Phi_{\sE}$
equals~$\Phi_{(a,\alpha)}$, which is defined by~\thetag{2.1}.

\proclaim{Proposition 2.14}
Suppose that $A=B$. We define the equivalence $\Phi_{(a,\alpha)}$
by \thetag{2.1} and consider the object $S_{(a,\alpha)}$ on
$A\myt\nomathbreak A$ that represents $\Phi_{(a,\alpha)}$. Then
the sheaf $\sJ(S_{(a,\alpha)})$ is equal to $\Delta_*P_{(\alpha,-a)}$,
where $\Delta$ is the diagonal embedding of $A\myt\hat A$ into
$(A\myt\hat A)\myt(A\myt\hat A)$ and $P_{(\alpha,-a)}$ is the line bundle
on $A\myt\hat A$ from Definition~$2.1$.
\endproclaim

\demo{Proof} By Assertion~2.8, the functor $\Phi_{S_A}$ maps the structure
sheaf $\sO_{(a',\alpha')}$ to the sheaf~$S_{(a',\alpha')}$ on $A\myt A$
(see~\thetag{2.2}). The functor~$\Ad_{S_{(a,\alpha)}}$ maps the sheaf
$S_{(a',\alpha')}$ to itself. Indeed, \thetag{1.6} implies that the object
$\Ad_{S_{(a,\alpha)}}(S_{(a',\alpha')})$ represents the functor
$$
\Phi_{(a,\alpha)}\circ\Phi_{(a',\alpha')}\circ\Phi^{-1}_{(a,\alpha)},
$$
which is in turn isomorphic to $\Phi_{(a',\alpha')}$ since all such functors
commute. We thus see that the functor determined by the sheaf
$\sJ(S_{(a,\alpha)})$ maps the structure sheaf of any point to itself,
whence it is a line bundle~$L$ supported on the diagonal.

To find the line bundle $L$, we find the image of the bundle
$P_{(\alpha',a')}$ under this functor. Using Assertion~2.8, we see that
$\Phi_{S_A}$ maps the bundle~$P_{(\alpha',a')}$ to the object
$\sO_{\{-a'\}\myt A}\ot p_2^*(P_{\alpha'})[n]$. It is easily verified that
the further action of~$\Ad_{S_{(a,\alpha)}}$ maps this object to the
object $\sO_{\{-a'+a\}\myt A}\ot p_2^*(P_{\alpha'+\alpha})[n]$. Hence
the functor determined by the sheaf $\sJ(S_{(a,\alpha)})$ maps the bundle
$P_{(\alpha',a')}$ to the bundle $P_{(\alpha'+\alpha,a'-a)}$. Therefore
the bundle $L$ is isomorphic to~$P_{(\alpha, -a)}$.
\enddemo

Given abelian varieties $A$ and $B$, we denote by $\sE q(A,B)$ the set of
all exact equivalence from the category~$\db A$ to the category~$\db B$ up
to an isomorphism.

We consider two groupoids $\fA$ and $\fD$ (that is, categories all of whose
morphisms are invertible). The objects of both categories are abelian
varieties. The morphisms in~$\fA$ are isomorphisms between the abelian
varieties as algebraic groups. The morphisms in $\fD$ are exact
equivalences between the derived categories of coherent sheaves on abelian
varieties, that is,
$$
\align
\lM or_{\fA}(A,B)
&:= \Iso(A,B),
\\
\lM or_{\fD}(A,B)
&:=\sE q(A,B).
\endalign
$$

By Theorem 2.10 there is a map of the set $\sE q(A,B)$ to the set
$\Iso(A\myt\hat A,B\myt\wht B)$, which sends each equivalence $\Phi_{\sE}$
to the isomorphism~$f_{\sE}$. We consider the map~$F$ from~$\fD$ to~$\fA$
that sends an abelian variety $A$ to the variety $A\myt\hat A$ and acts
on the morphisms as described above.

\proclaim{Proposition 2.15}
The map $F\:\fD\to\fA$ is a functor.
\endproclaim

\demo{Proof}
To prove this, we must verify that $F$ preserves composition of morphisms.
Consider abelian varieties $A$, $B$, $C$. Let $\sE$ and $\sF$ be objects
of the categories $\db{A\myt B}$ and $\db{B\myt C}$ respectively such that
the functors
$$
\align
&\Phi_{\sE}\:\db A\to\db B,
\\
&\Phi_{\sF}\:\db B\to\db C
\endalign
$$
are equivalences. We denote by $\sG$ the object in $\db{A\myt C}$ that
represents the composite of these functors.

The relation~\thetag{1.4} yields an isomorphism
$\Ad_{\sG}\cong\Ad_{\sF}\circ\Ad_{\sE}$, whence we get
$$
\alignat1
\Phi_{\!\sJ(\sF)}\circ\Phi_{\!\sJ(\sE)}
&\cong(\Phi_{S_A}^{-1}\circ\Ad_{\sF}\circ\Phi_{S_A})
\circ(\Phi_{S_A}^{-1}\circ\Ad_{\sE}\circ\Phi_{S_A})
\\
&\cong\Phi_{S_A}^{-1}\circ\Ad_{\sG}\circ\Phi_{S_A}\cong\Phi_{\!\sJ(\sG)}.
\endalignat
$$
By Theorem 2.10, the objects $\sJ(\sE)$, $\sJ(\sF)$, $\sJ(\sG)$ are line
bundles supported on the graphs of the isomorphisms $f_{\sE}$,
$f_{\sF}$, $f_{\sG}$. We thus get the equation
$f_{\sG}=f_{\sF}\cdot f_{\sE}$.
\enddemo

\proclaim{Corollary 2.16}
Suppose that $A$ is an abelian variety and $\Phi_{\sE}$ is an
autoequivalence of the derived category~$\db A$.
Then the map $\Phi_{\sE}\mapsto f_{\sE}$ determines a group homomorphism
$$
\gamma_A\:\Autoeq\db A\to\Aut(A\myt\hat A).
$$
\endproclaim

We thus have the functor $F\:\fD\to\fA$. Our next purpose is to describe
it. To do this, we must study which elements of
$\Iso(A\myt\hat A,B\myt\wht B)$ may be realized as~$f_{\sE}$ for some
$\sE$ and which equivalences~$\sE_1$ and~$\sE_2$ satisfy the equation
$f_{\sE_1}=f_{\sE_2}$.

We consider an arbitrary morphism $f\:A\myt\hat A\to B\myt\wht B$.
It is convenient to write it as a matrix
$$
\pmatrix
\alpha & \beta
\\
\gamma & \delta
\endpmatrix,
$$
where the morphism $\alpha$ maps $A$ to $B$, \ $\beta$ maps $\hat A$ to
$B$, \ $\gamma$ maps $A$ to~$\wht B$, and $\delta$ maps $\hat A$ to~$\wht B$.
Each morphism~$f$ determines two other morphisms
$\hat f$ and~$\myw f$ from $B\myt\wht B$ to $A\myt\hat A$
whose matrices are
$$
\hat f=\pmatrix
\wht\delta & \wht\beta
\\
\wht\gamma & \wht\alpha
\endpmatrix, \qquad
\myw f=\pmatrix
\wht\delta & -\wht\beta
\\
-\wht\gamma & \wht\alpha
\endpmatrix.
$$

We define a set $U(A\myt\hat A,B\myt\wht B)$ as the subset of all $f$
in $\Iso(A\myt\hat A,B\myt\wht B)$ such that $\myw f$ coincides with the
inverse to~$f$, that is,
$$
U(A\myt\hat A,B\myt\wht B)
:=\{f\in\Iso(A\myt\hat A,B\myt\wht B)|\myw f=f^{-1}\}.
$$
If $B=A$, we denote this set by $U(A\myt\hat A)$. We note that
$U(A\myt\hat A)$ is a subgroup in $\Aut(A\myt\hat A)$.

\definition{Definition 2.17}
An isomorphism $f\:A\myt\hat A\mov\sim\to\to B\myt\wht B$ is called {\it
isometric} if it belongs to $U(A\myt\hat A,B\myt\wht B)$.
\enddefinition

\proclaim{Proposition 2.18}
For every equivalence $\Phi_{\sE}\:\db A\mov\sim\to\to\db B$,
the isomorphism $f_{\sE}$ is isometric.
\endproclaim

\demo{Proof}
Passing to the algebraic closure if necessary, we may assume that
the field $k$ is algebraically closed. To verify the equation
$\myw f_{\sE}=f^{-1}_{\sE}$, it suffices to prove that these morphisms
coincide at closed points. Suppose that~$f_{\sE}$ sends a point
$(a,\alpha)\in A\myt\hat A$ to a point $(b,\beta)\in B\myt\wht B$.
We claim that $\myw f_{\sE}(b,\beta)=(a,\alpha)$ or, equivalently, that
$\hat f(-b,\beta)=(-a,\alpha)$.

The isomorphism $f_{\sE}$ is determined by an abelian subvariety
$X\hookrightarrow A\myt\hat A\myt B\myt\wht B$. Hence we must verify that
$P_{(0,0,\beta,-b)}\ot\sO_X\cong P_{(\alpha,-a,0,0)}\ot\sO_X$ or,
equivalently, that the sheaf $\sJ':=P_{(-\alpha,a,\beta,-b)}\ot\sJ(\sE)$
is isomorphic to the sheaf~$\sJ(\sE)$.

By Proposition 2.14, the functor determined by $\sJ'$ is a composite of
the functors represented by the objects
$\sJ(S_{(-a,-\alpha)})$, $\sJ(\sE)$ and
$\sJ(S_{(b,\beta)})$. Thus $\sJ'$ coincides with~$\sJ(\sE')$,
where $\sE'$ is the object in $\db{A\myt B}$ that represents the functor
$$
\Phi_{(b,\beta)}\circ\Phi_{\sE}\circ\Phi_{(-a,-\alpha)}.
$$
This composite is isomorphic to the functor $\Phi_{\sE}$ by Corollary~2.13.
Hence the object $\sE'$ is isomorphic to~$\sE$,
and we have $\sJ'=\sJ(\sE')\cong\sJ(\sE)$.
\enddemo

As a corollary of Theorem 2.10 and Proposition 2.18,
we get the following result.

\proclaim{Theorem 2.19}
Let $A$, $B$ be abelian varieties over a field $k$. If the derived
categories of coherent sheaves $\db A$ and~$\db B$ are equivalent as
triangulated categories, then there is an isometric isomorphism between
$A\myt\hat A$ and $B\myt\wht B$.
\endproclaim

The converse assertion is also true for abelian varieties over an
algebraically closed field of characteristic $0$~\cite{13}. Another
proof of this result is given in~\S\,4.

\proclaim{Corollary 2.20}
For every abelian variety $A$ there are only finitely many non-isomorphic
abelian varieties whose derived categories of coherent sheaves are
equivalent to $\db A$ (as triangulated categories).
\endproclaim

\demo{Proof}
As proved in~\cite{6}, for every abelian variety $Z$ there are only finitely
many (up to an isomorphism) abelian varieties that can be embedded into~$Z$
as an abelian subvariety. Applying this assertion to $Z=A\myt\hat A$ and
using Theorem~2.19, we get the required result.
\enddemo

In conclusion we would like to explain the term ``isometric" for
elements of $U(A\myt\hat A,B\myt\wht B)$. Suppose that $k$
is the field of complex numbers. We denote the first homology lattices
$H_1(A,\ZZ)$ and $H_1(B,\ZZ)$ by $\varGamma_A$ and $\varGamma_B$
respectively. Any lattice representable as $\varGamma\oplus\varGamma^*$
carries a canonical symmetric bilinear form
$$
Q\bigl((x,l),(y,m)\bigr)=l(y)+m(x).
$$
We denote by $Q_A$ and $Q_B$ the corresponding symmetric bilinear forms on
$\varLambda_A:=H_1(A\myt\hat A,\ZZ)$ and $\varLambda_B:=H_1(B\myt\wht B,\ZZ)$.

The set of homomorphisms from $A\myt\hat A$ to $B\myt\wht B$ is a subset
of $\Hom_{\ZZ}(\varLambda_A,\varLambda_B)$. The elements of
$U(A\myt\hat A,B\myt\wht B)$ can be described in these terms as follows.

\proclaim{Proposition 2.21}
An isomorphism $f\:A\myt\hat A\to B\myt\wht B$ belongs to
$U(A\myt\hat A,B\myt\wht B)$ if and only if it determines an isometry
of the lattices $(\varLambda_A,Q_A)$ and $(\varLambda_B,Q_B)$,
that is, $F^tQ_BF=Q_A$, where $F\:\varLambda_A\to\varLambda_B$
is the map induced by $f$ on the first homology.
\endproclaim

\demo\nofrills{The proof} \ is obtained by a direct matrix computation
using Remark~2.3.
\enddemo

\head
\S\,3. Objects that represent equivalences, and the autoequivalence groups
\endhead

Propositions 2.15,~2.18 imply that there is a homomorphism from the group
$\Auteq\db A$ of exact autoequivalences to the group $U(A\myt\hat A)$ of
isometric isomorphisms. In this section we describe the kernel of this
homomorphism. By Proposition~2.14 we know that all equivalences
$\Phi_{(a,\alpha)}[n]$ (see~\thetag{2.1}) belong to the kernel.
We shall show that they exhaust the kernel. To prove this, we use
an assertion of independent interest: for abelian varieties, if the
functor~$\Phi_{\sE}$ is an equivalence, then the object $\sE$ is a
sheaf on the product (up to a shift in the derived category).
This assertion is specific to abelian varieties and breaks in other
cases, say, for K3-surfaces.

\proclaim{Lemma 3.1}
Let $\sE$ be an object on the product $A\myt B$ such that
$\Phi_{\sE}\:\db A\to\db B$ is an equivalence. We consider the projection
$q\:(A\myt\hat A)\myt(B\myt\nomathbreak\wht B)\to A\myt B$
and denote by~$K$ the direct image $\bR^\cdot q_*\sJ(\sE)$,
where~$\sJ(\sE)$ is the object from Definition~$2.9$.
Then $K$ is isomorphic to the object $\sE\ot(\sE^\vee|_{(0,0)})$,
where $\sE^\vee|_{(0,0)}$ denotes the complex of vector spaces
which is the inverse image of the object
$\bR^\cdot\HHom(\sE,\sO_{A\myt B})$ under the embedding of the point
$(0,0)$ to the abelian variety $A\myt B$.
\endproclaim

\demo{Proof} We consider the abelian variety
$$
Z=(A\myt\hat A)\myt(A\myt A)\myt(B\myt B)\myt(B\myt\wht B)
$$
and the object
$$
H=p_{1234}^*S_A\otimes p^*_{35}\sE^\vee[n]\otimes p^*_{46}\sE
\otimes p_{5678}^*S_B^\vee[2n].
$$
Using Proposition 1.1 on the composition of functors and
diagram~\thetag{2.3}, we see that $\sJ(\sE)\cong p_{1278*}H$.
Hence, the object~$K$ equals~$p_{17*}H$. To compute the last object,
we consider the projection of~$Z$ onto
$$
V=A\myt(A\myt A)\myt(B\myt B)\myt B
$$
and denote it by~$v$. To calculate $v_*H$, we use the fact that the functor
$\Phi_{S_A}$ is the composite of $\Phi_{\sP_A}$ and $\bR\mu_{A*}$, where
$$
\sP_A=p^*_{14}\sO_\Delta\ot p^*_{23}P\in\db{(A\myt\hat A)\myt(A\myt A)}.
$$
It is easy to see that $p_{134*}\sP_A\cong\sO_{T_A}[-n]$, where
$T\subset A\myt A\myt A$ is the subvariety isomorphic to~$A$ and consisting
of the points $(a,0,a)$. Taking into account that
$\mu_A(a_1,a_2)=\bigl(a_1,m(a_1,a_2)\bigr)$, we see that $p_{134*}S_A$
is also isomorphic to~$\sO_{T_A}[-n]$. We similarly verify that
$p_{134*}S_B^\vee[2n]=\sO_{T_B}$.

Thus we have
$$
v_*H\cong p_{123}^*\sO_{T_A}\otimes p^*_{24}\sE^\vee\otimes p^*_{35}\sE
\otimes p_{456}^*\sO_{T_B}
$$
on $V$. We consider the embedding
$$
j\:A\myt A\myt B\myt B\to V, \qquad
(a_1,a_2,b_1,b_2)\mapsto(a_1,0,a_2,0,b_1,b_2).
$$
The object $v_* H$ is isomorphic to $j_* M$, where
$$
M=(\sE^\vee|_{(0,0)})\otimes p^*_{12}\sO_{\Delta_A}\otimes p_{23}^*\sE
\otimes p_{34}^*\sO_{\Delta_B}.
$$
We finally obtain that
$K\cong p_{14*}M\cong(\sE^\vee|_{(0,0)})\otimes\sE$.
\enddemo

\proclaim{Proposition 3.2}
Let $A$, $B$ be abelian varieties, and let $\sE$ be an object of
$\db{A\myt B}$ such that the functor $\Phi_{\sE}\:\db A\mov\sim\to\to\db B$
is an exact equivalence. Then $\sE$ has only one non-trivial cohomology,
that is, $\sE$ is isomorphic to the object~$\sF[n]$, where~$\sF$ is a sheaf
on $A\myt B$.
\endproclaim

\demo{Proof} Consider the projection
$$
q\:(A\myt\hat A)\myt(B\myt\wht B)\to A\myt B
$$
and denote by $q'$ its restriction to the abelian subvariety~$X$ which is
the support of the sheaf~$\sJ(\sE)$ and the graph of the isomorphism
$f_{\sE}$. By Theorem~2.10, the sheaf~$\sJ(\sE)$ equals~$i_*(L)$, where
$L$ is a line bundle on~$X$.

We denote the object $\bR^\cdot q_*\!\sJ(\sE)=\bR^\cdot q'_*L$ by~$K$.
The morphism $q'$ is a homomorphism of abelian varieties. Let $d$ be the
dimension of~$\Ker (q')$. Then $\dim\Im(q')=2n-d$, whence the cohomology
sheaves~$H^j(K)$ are trivial for $j\notin[0,d]$.

On the other hand, the object $K$ is isomorphic to
$\sE\ot(\sE^\vee|_{(0,0)})$ by Lemma~3.1.

Shifting the object $\sE$ in the derived category if necessary, we may
assume that the highest non-zero cohomology of~$\sE$ is~$H^0 ({\sE})$.
Let $H^{-i}({\sE})$ be the lowest non-zero cohomology of $\sE$
(here $i\ge0$), and let $H^k({\sE}^\vee)$ be the highest non-zero
cohomology of~$\sE^\vee$. Replacing $\sE$ by~$T^*_{(a,b)}{\sE}$
if necessary, we may assume that the point $(0,0)$ belongs to the support
of the sheaf~$H^k({\sE}^\vee)$. Since the supports of $\sE$ and~$K$
coincide, the supports of all cohomology sheaves of~$\sE$ belong to
$\Im(q')$. In particular, we have $\codim\Supp H^{-i}({\sE})\ge d$.
Hence the object $\bigl(H^{-i}({\sE})\bigr)^\vee[-i]$
has trivial cohomologies in all degrees less than $i+d$.

The canonical morphism $H^{-i}({\sE})[i]\to\sE$ induces a non-trivial
morphism
$$
\sE^\vee\to\bigl(H^{-i}({\sE})\bigr)^\vee[-i].
$$
Since the degrees of non-trivial cohomologies of the second object belong
to $[i+d,\infty)$, we see that $k\ge i+d$, where~$H^k({\sE}^\vee)$ is the
highest non-zero cohomology of~$\sE^\vee$. Thus we obtain that the object
$$
K=\sE^\vee|_{(0,0)}\ot\sE
\tag3.1
$$
has a non-trivial cohomology at the same degree $k\ge i+d$. On the other
hand, it is known that the cohomology sheaves $H^j(K)$ are trivial for
$j\notin[0,d]$. This is possible only when $i=0$. Hence the object~$\sE$
has only one non-trivial cohomology at degree~$0$ and, therefore, $\sE$
is isomorphic to a sheaf.
\enddemo

We consider the case $B\cong A$. Let $\sE$ be a sheaf on $A\myt A$ such that
$\Phi_{\sE}$ is an autoequivalence. Let us describe all $\sE$ such that
$f_{\sE}$ is the identity map and, therefore, the graph~$X$ of this map is
the diagonal in $(A\myt\hat A)\myt(A\myt\hat A)$. We thus obtain that the
object
$$
K=\sE^\vee|_{(0,0)}\ot\sE=\bR^\cdot q_*\sJ(\sE)
$$
takes the form $\Delta_*(M)$, where $M$ is an object on $A$ and
$\Delta\: A\to A\myt A$ is the diagonal embedding.

We suppose that the point $(0,0)$ belongs to the support of the sheaf ${\sE}$.
Hence $\sE^\vee|_{(0,0)}$ is a non-trivial complex of vector spaces. Then
the assumption $K=\Delta_*(M)$ implies that there is a sheaf~$E$ on~$A$
such that $\sE\cong\Delta_*(E)$. Therefore
$\Phi_{\sE}(\cdot)\cong E\ot(\cdot)$. Since~$\Phi_{\sE}$ is an
autoequivalence, $E$ is a line bundle. We easily verify that
the condition $f_{\sE}=\id$ may hold only when $E\in\Pic^0(A)$.

If the point $(0,0)$ does not belong to $\Supp\sE$, then we can replace
$\sE$ by a sheaf $\sE':=T_{(a_1,a_2)*}\sE$ whose support contains $(0,0)$.
Proposition~2.14 implies that $f_{\sE'}=f_{\sE}$. As shown above, there is
an isomorphism $\sE'\cong\Delta_*(E')$, where $E'\in\Pic^0(A)$.
Hence $\sE\cong T_{(a_1-a_2,0)*}\Delta_*(E')$. Thus we get the following
proposition.

\proclaim{Proposition 3.3}
Let $A$ be an abelian variety. Then the kernel of the homomorphism
$$
\gamma_A\:\Auteq\db A\to U(A\myt\hat A)
$$
consists of the autoequivalences $\Phi_{(a,\alpha)}[i]=
T_{a*}(\cdot)\ot P_\alpha[i]$ and is thus isomorphic to the group
$\ZZ\oplus(A\myt\hat A)_k$, where $(A\myt\hat A)_k$ is the group of
$k$-points of the abelian variety $A\myt\hat A$.
\endproclaim

\proclaim{Corollary 3.4}
Let $A$, $B$ be abelian varieties, and let
$\sE_1$,~$\sE_2$ be objects on the product
$A\myt B$ that determine equivalences between the derived categories of
coherent sheaves. If $f_{\sE_1}=f_{\sE_2}$, then
$$
\sE_2\cong T_{a*}\sE_1\otimes P_\alpha[i]
$$
for some $k$-point $(a,\alpha)\in A\myt\hat A$.
\endproclaim

\head
\S\,4. Semihomogeneous vector bundles
\endhead

As shown in previous sections, an equivalence $\Phi_{\sE}$ from~$\db A$
to~$\db B$ induces an isometric isomorphism of the varieties
$A\myt\hat A$ and $B\myt\wht B$. In this section we suppose that
the field~$k$ is algebraically closed and $\har(k)=\nomathbreak0$.
Using the technique of~\cite{10} and the results of~\cite{1},
we shall show under these assumptions that every isometric isomorphism
$f\:A\myt\hat A\to B\myt\wht B$ is induced by some equivalence $\Phi_{\sE}$.
The fact that the existence of an isometric isomorphism between
$A\myt\hat A$ and $B\myt\wht B$ implies equivalence of the derived
categories~$\db A$ and~$\db B$ was proved in~\cite{13}. We shall
give another proof of this fact.

We recall that any line bundle $L$ on an abelian variety~$D$ determines a
map~$\phi_L$ from~$D$ to~$\wht D$, which sends a point~$d$ to the point
corresponding to the bundle $T^*_dL\ot L^{-1}\in\Pic^0(D)$. This
correspondence yields an embedding of~$\NS(D)$ into $\Hom(D,\wht D)$.
Moreover, it is known that a map $\phi\:D\to\wht D$ belongs to the image
of~$\NS(D)$ if and only if $\wht\phi=\phi$.

Semihomogeneous bundles on abelian varieties enable us to generalize this
phenomenon as follows. Every element of $\NS(D)\ot\QQ$ determines some
correspondence on $D\myt\wht D$, and every such correspondence is obtained
from some semihomogeneous bundle (see Proposition~4.6 below and Lemma~2.13).
We recall some results on homogeneous and semihomogeneous bundles on
abelian varieties.

\definition{Definition 4.1}
A vector bundle $E$ on an abelian variety~$D$ is called {\it homogeneous}
if $T^*_d(E)\cong E$ for every point $d\in D$.
\enddefinition

\definition{Definition 4.2}
A vector bundle $F$ on an abelian variety~$D$ is called {\it unipotent}
if there is a filtration
$$
0=F_0\subset F_1\subset\cdots\subset F_n=F
$$
such that $F_i/F_{i-1}\cong\sO_D$ for all $i=1,\dots,n$.
\enddefinition

The following proposition characterizes homogeneous vector bundles.

\proclaim{Proposition 4.3~\cite{8},~\cite{10}}
Let $E$ be a vector bundle on an abelian variety~$D$. Then the following
conditions are equivalent.

{\rm(i)} The vector bundle $E$ is homogeneous.

{\rm(ii)} There are line bundles $P_i\in\Pic^0(D)$ and unipotent bundles
$F_i$ such that $E\cong\bigoplus_i(F_i\ot P_i)$.
\endproclaim

\definition{Definition 4.4}
A vector bundle $E$ on an abelian variety~$D$ is called {\it
semihomogeneous} if, for every point $d\in D$, there is a line bundle
$L$ on~$D$ such that $T^*_d(E)\cong E\ot L$. (We note that then~$L$
belongs to $\Pic^0(D)$.)
\enddefinition

We recall that a vector bundle on a variety is called
{\it simple} if its automorphism algebra coincides with the field~$k$.

\proclaim{Proposition 4.5 (\cite{10}, Theorem~5.8)}
Let $E$ be a simple vector bundle on an abelian variety~$D$.
Then the following conditions are equivalent.

{\rm1)} $\dim H^j\bigl(D,\End(E)\bigr)=\binom nj$ for every $j=0,\dots,n$.

{\rm2)} $E$ is a semihomogeneous bundle.

{\rm3)} $\End(E)$ is a homogeneous bundle.

{\rm4)} $E\cong\pi_*(L)$ for some isogeny $\pi\:Y\to D$ and some
line bundle $L$ on $Y$.
\endproclaim

Let $E$ be a vector bundle on an abelian variety $D$. We denote by~$\mu(E)$
the equivalence class of $\frac{\det(E)}{r(E)}$ in
$\NS(D)\otimes_\ZZ\QQ$.
Every element $\mu=\frac{[L]}l\in\NS(D)\otimes_\ZZ\QQ$
(hence every bundle~$E$) determines a correspondence
$\Phi_\mu\subset D\myt\wht D$ by
$$
\Phi_\mu=\Im\bigl[D\mov{(l,\phi_L)}\to\to D\myt\wht D\bigr].
$$
Here $\phi_L$ is the well-known map from $D$ to $\wht D$ which sends a
point~$d$ to the point corresponding to the bundle
$T^*_dL\ot L^{-1}\in\Pic^0(D)$. We denote the projections of
$\Phi_\mu$ onto $D$ and $\wht D$ by $q_1$ and $q_2$ respectively.
In the particular case when the bundle is a line bundle~$L$, we get the
graph of the map $\phi_L\:D\to\wht D$.

A complete description of all simple semihomogeneous bundles is given
in~\cite{10}.

\proclaim{Proposition 4.6 (\cite{10}, Theorem 7.10)}
Suppose that $\mu=\frac{[L]}l$, where $[L]$ is the equivalence class of
the bundle~$L$ in~$\NS(D)$ and $l$ is a positive integer. Then

{\rm1)} there is a simple semihomogeneous vector bundle $E$ of slope
$\mu(E)=\mu$;

{\rm2)} every simple semihomogeneous vector bundle $E'$ of slope
$\mu(E')=\mu$ coincides with $E\ot M$ for some line bundle
$M \in\Pic^0(D)$;

{\rm3)} we have equations $r(E)^2=\deg(q_1)$ and $\chi(E)^2=\deg(q_2)$.
\endproclaim

The following assertion enables us to characterize all semihomogeneous
vector bundles in terms of simple bundles.

\proclaim{Proposition 4.7 (\cite{10}, Propositions 6.15, 6.16)}
Every semihomogeneous vector bundle $F$ of slope~$\mu$ admits a filtration
$$
0=F_0\subset F_1\subset\cdots\subset F_n=F
$$
such that $E_i=F_i/F_{i-1}$ are simple semihomogeneous vector bundles
of the same slope~$\mu$. Every simple semihomogeneous bundle is stable.
\endproclaim

We shall use the following two lemmas on semihomogeneous bundles. These
lemmas are direct corollaries of the previous assertions.

\proclaim{Lemma 4.8}
Any two simple semihomogeneous bundles $E_1$, $E_2$ of the same slope~$\mu$
are either isomorphic or orthogonal to each other, that is, either
$E_1=E_2$ or
$$
\Ext^i(E_1,E_2)=0, \qquad
\Ext^i(E_2,E_1)=0
$$
for all $i$.
\endproclaim

\demo{Proof}
Proposition 4.6 implies that $E_2\cong E_1\ot M$, whence $\HHom(E_1,E_2)$
is a homogeneous bundle. By Proposition~4.3, every homogeneous bundle is
equal to a direct sum of unipotent bundles twisted by line bundles
belonging to $\Pic^0(D)$. Therefore, either all cohomologies of
$\HHom(E_1,E_2)$ are equal to zero (and hence the bundles $E_1$,~$E_2$
are orthogonal), or the bundle $\HHom(E_1,E_2)$ has a non-zero section.
In the last case we get a non-zero homomorphism from
$E_1$ to~$E_2$. But these two bundles are stable and have the same slope.
Hence every non-zero homomorphism is actually an isomorphism.
\enddemo

\proclaim{Lemma 4.9}
Let $E$ be a simple semihomogeneous vector bundle on an abelian variety~$D$.
We have $T^*_d(E)\cong E\ot P_\delta$ if and only if
$(d,\delta )\in\Phi_\mu$.
\endproclaim

\demo{Proof}
Let us prove that $T^*_d(E)\cong E\ot P_\delta$ for every point
$(d,\delta)\in\Phi_\mu$. Indeed, we put $l=r(E)$ and $L=\det(E)$.
By definition of $\Phi_\mu$ we can write
$(d,\delta)=\bigl(lx,\phi_L(x)\bigr)$ for some point $x\in D$.
Since $E$ is semihomogeneous, there is a line bundle
$M\in\Pic^0(D)$ such that
$$
T^*_x (E)\cong E\ot M.
\tag 4.1
$$
Comparing the determinants, we get $T^*_x(L)\cong L\ot M^{\ot l}$.
By definition of the map~$\phi_L$, this means that $P_{\phi_L(x)}=M^{\ot l}$.
On the other hand, iterating the equation~\thetag{4.1} $l$ times, we get
$$
T^*_{lx}(E)\cong E\ot M^{\ot l}=E\ot P_{\phi_L(x)}.
$$
Therefore $T^*_d(E)\cong E\ot P_\delta$ because
$(d,\delta)=\bigl(lx,\phi_L(x)\bigr)$.

Conversely, define a subgroup $\Sigma^0(E)\subset\wht D$ by the condition
$$
\Sigma^0(E):=\{\delta\in\wht D\mid E\ot P_\delta\cong E\}.
\tag4.2
$$
Since $E$ is semihomogeneous, the bundle $\End(E)$ is
homogeneous by Proposition~4.5.
Thus $\End(E)$ may be presented as a sum $\bigoplus_i(F_i\ot P_i)$,
where all~$F_i$ are unipotent. Hence
$H^0\bigl(\End(E)\ot P\bigr)\ne0$ for at most~$r^2$ line bundles
$P\in\Pic^0(D)$, that is, the order of the group $\Sigma^0(E)$ does not
exceed~$r^2$. On the other hand, it is known that
$q_2\bigl(\Ker(q_1)\bigr)\subset\Sigma^0(E)$.
Hence we obtain that $\ord\Sigma^0(E)=r^2$ and
$q_2\bigl(\Ker(q_1)\bigr)=\Sigma^0(E)$.

We now suppose that $T^*_d(E)\cong E\ot P_\delta$ for some point
$(d,\delta )\in D\myt\wht D$. Consider the point $\delta'\in\wht D$
such that $(d,\delta')\in\Phi_\mu$. It was shown that then we have
an isomorphism $T^*_d(E)\cong E\ot P_{\delta'}$. Hence
$E\ot P_{(\delta-\delta')}\cong E$ and, therefore,
$(\delta-\delta')\in\Sigma^0(E)$.
Since $\Sigma^0(E)=q_2\bigl(\Ker(q_1)\bigr)$, the point
$(0,\delta-\delta')$ belongs to~$\Phi_\mu$. Hence the point
$(d,\delta)$ also belongs to~$\Phi_\mu$.
\enddemo

We now present a construction that starts from an isometric isomorphism
$f$ and produces an object~$\sE$ on the product such that $\sE$ determines
an equivalence of the derived categories and $f_{\sE}$ coincides with~$f$.

\example{Construction 4.10}
We fix an isometric isomorphism $f\:A\myt\hat A\to B\myt\wht B$ and denote
its graph by~$\Gamma$. As above, we write~$f$ in the matrix form:
$$
f=\pmatrix
x & y
\\
z & w
\endpmatrix.
$$
Suppose that $y\:\hat A\to B$ is an isogeny. Then the map~$f$ determines
an element $g\in\Hom(A\myt B,\hat A\myt\wht B)\ot_{\ZZ}\QQ$ by the formula
$$
g=\pmatrix
y^{-1}x & -y^{-1}
\\
-\hat y^{-1} & wy^{-1}
\endpmatrix.
$$

The element $g$ determines a correspondence on
$(A\myt B)\myt(\hat A\myt\wht B)$.
Since $f$ is isometric, we easily verify that $\hat g=g$. This means that
$g$ actually belongs to the image of
$\NS(A\myt B)\ot_{\ZZ}\QQ$ under the canonical embedding into
$\Hom(A\myt B,\hat A\myt\wht B)\ot_\ZZ\QQ$ (see, for example,~\cite{11}).
Hence there is $\mu=\frac{[L]}l\in\NS(A\myt B)$ such that $\Phi_\mu$
coincides with the graph of the correspondence~$g$. Proposition~4.6
asserts that for every~$\mu$ one can construct a simple semihomogeneous
bundle~$\sE$ on $A\myt B$ of slope $\mu(\sE)=\mu$.

Below, we shall see that the functor $\Phi_{\sE}$ from $\db A$ to $\db B$
is an equivalence and $f_\sE=f$. But now we compare the graphs
$\Gamma$ and~$\Phi_\mu$.
If a point $(a,\alpha,b,\beta)$ belongs to~$\Gamma$, then
$$
b=x(a)+y(\alpha), \qquad
\beta=z(a)+w(\alpha)
$$
and, therefore,
$$
\align
\alpha
&=-y^{-1}x(a)+y^{-1}(b),
\\
\beta
&=(z- wy^{-1}x)(a)+wy^{-1}(b).
\endalign
$$

Since $f$ is isometric, we have $(z-wy^{-1}x)=-\hat y^{-1}$. Hence the point
$(a,\alpha,b,\beta)$ belongs to $\Gamma$ if and only if
$(a,-\alpha,b,\beta)$ belongs to~$\Phi_\mu$. We thus find that
$\Phi_\mu=(1_A,-1_{\hat A},1_B,1_{\wht B})\Gamma$.

In particular, since $f$ is an isomorphism, the projections of~$\Phi_\mu$
onto $A\myt\hat A$ and $B\myt\wht B$ are also isomorphisms.
\endexample

\proclaim{Proposition 4.11}
Let $\sE$ be the semihomogeneous bundle on $A\myt B$ constructed from an
isometric isomorphism~$f$ in the way described above. Then the functor
$\Phi_{\sE}\:\db A\to\db B$ is an equivalence.
\endproclaim

\demo{Proof}
We denote by $\sE_a$ the restriction of the bundle $\sE$ to the fibre
$\{a\}\myt B$. By Theorem~1.5, to prove that~$\Phi_{\sE}$ is fully
faithful, we must verify that all the bundles~$\sE_a$ are simple and
orthogonal to each other for different points.

We first note that, by Proposition~4.6, the rank of the bundle~$\sE$
equals the square root of the degree of the map $\Phi_\mu\to A\myt B$,
that is, $\sqrt{\deg(\beta)}$.

It follows from the semihomogeneity of $\sE$ that all bundles~$\sE_a$
are semihomogeneous. Moreover, the slope $\mu(\sE_a)$ of the restriction
is equal to $\delta\beta^{-1}\in\NS(B)\ot\nomathbreak\QQ
\subset\Hom(B,\wht B)\ot\QQ$. For brevity, we denote the element
$\delta\beta^{-1}$ by $\nu$ and regard it as an element of $\NS(B)\ot\QQ$.
By Proposition~4.6 there is a simple semihomogeneous bundle $F$ on~$B$
of slope $\mu(F)=\nu$. Then $\Phi_\nu$ is clearly equal to
$\Im\bigl[\hat A\mov{(\beta,\delta)}\to\to B\myt\wht B\bigr]$.
Since $f$ is an isomorphism, the map
$\hat A\mov{(\beta,\delta)}\to\to B\myt\wht B$ is an embedding.
Using Proposition~4.6 again, we get the equations
$r(F)=\sqrt{\deg(\beta)}=r(\sE_a)$. Thus the bundles $F$ and $\sE_a$
are semihomogeneous, of the same slope and rank. Moreover, $F$ is simple.
It follows form Propositions~4.7 and~4.6,~2) that the bundle $\sE_a$ is
also simple.

Lemma 4.8 now implies that for any points $a_1,a_2\in A$ the bundles
$\sE_{a_1}$ and $\sE_{a_2}$ are either orthogonal or isomorphic.
Suppose that they are isomorphic. Since $\sE$ is semihomogeneous, we have
$$
T^*_{(a_2-a_1,0)}\sE\cong\sE\ot P_{(\alpha,\beta)}
\tag4.3
$$
for some point $(\alpha,\beta)\in\hat A\myt\wht B$. In particular,
$$
\sE_{a_2}\ot P_\beta\cong\sE_{a_1}\cong\sE_{a_2}.
$$
Therefore $P_\beta\in\Sigma^0(\sE_a)$ (see~\thetag{4.2}).

By Lemma 4.9 and Proposition 4.6, the orders of $\Sigma^0(E)$ and
$\Sigma^0(E_a)$ are equal to~$r^2$. We claim that the natural map
$\sigma\:\Sigma^0(\sE)\to\Sigma^0(\sE_a)$ is an isomorphism. Indeed,
otherwise there is a point $\alpha'\in\hat A$ such that
$\sE\ot P_{\alpha'}\cong\sE$. Then Lemma~4.9 yields that
$(0,\alpha',0,0)\in\Phi_\mu$, contrary to the fact that the projection
$\Phi_\mu\to B\myt\wht B$ is an isomorphism.

Since $\sigma$ is an isomorphism, there is a point $\alpha'\in\hat A$
such that $\sE\ot P_{(\alpha',\beta)}\cong\sE$. It follows from
\thetag{4.3} that
$$
T^*_{(a_2-a_1,0)}\sE\cong\sE\ot P_{(\alpha-\alpha',0)}.
$$
By Lemma 4.9, this means that the point $(a_2-a_1,\alpha-\alpha',0,0)$
belongs to $\Phi_\mu$. Since the projection $\Phi_\mu\to B\myt\wht B$
is an isomorphism, we see that $a_2-a_1=0$. Thus the bundles $\sE_{a_1}$
and $\sE_{a_2}$ are orthogonal for different points $a_1$,~$a_2$.
Hence the functor $\Phi_{\sE}\:\db A\to\db B$ is fully faithful.
The same arguments show that the adjoint functor $\Psi_{\sE^\vee}$
is also fully faithful. Hence $\Phi_{\sE}$ is an equivalence.
\enddemo

\proclaim{Proposition 4.12}
Let $\sE$ be the semihomogeneous bundle constructed from an isometric
isomorphism $f\:A\myt\hat A\to B\myt\wht B$ in the  way described above.
Then $f_{\sE}=f$.
\endproclaim

\demo{Proof}
We denote by $X$ the graph of the morphism $f_{\sE}$. It follows from
Corollary 2.13 that a point $(a,\alpha,b,\beta)$ belongs to~$X$ if and only
if
$$
T_{b*}\sE\ot P_\beta\cong T^*_a\sE\ot P_\alpha,
$$
which is equivalent to the equation
$T^*_{(a,b)}\sE\cong E\ot P_{(-\alpha,\beta)}$.
Hence Lemma~4.9 implies that
$X=(1_A,-1_{\hat A},1_B,1_{\wht B})\Phi_\mu$, where $\mu=\mu(\sE)$ is the
slope of~$\sE$. On the other hand, the graph $\Gamma$ of the map~$f$ is
also equal to $(1_A,-1_{\hat A},1_B,1_{\wht B})\Phi_\mu$ by
Construction~4.10. Hence the isomorphisms $f_{\sE}$ and $f$ coincide.
\enddemo

In the construction of $\sE$ from $f$, we assumed that the
map $y\:\hat A\to B$ is an isogeny. If this is not the case, we
present~$f$ as a composite of two maps
$f_1\in U(A\myt\hat A,B\myt\wht B)$ and $f_2\in U(A\myt\hat A)$
such that $y_1$ and~$y_2$ are isogenies. We easily see that this can always
be done. For every~$f_i$, we find the corresponding object~$\sE_i$.
Then we consider the composite of the functors
$\Phi_{\sE_i}$ and take the object that represents it.

Assertions proved in this and previous sections can be joint into the
following theorems.

\proclaim{Theorem 4.13}
Let $A$, $B$ be abelian varieties over an algebraically closed field of
characteristic~$0$. Then the bounded derived categories of coherent sheaves
$\db A$ and $\db B$ are equivalent as triangulated categories if and only
if there is an isometric isomorphism
$f\:A\myt\hat A\mov\sim\to\to B\myt\wht B$.
\endproclaim

\proclaim{Theorem 4.14}
Let $A$ be an abelian variety over an algebraically closed field of
characteristic~$0$. Then the group $\Auteq\db A$ of exact autoequivalences
of the derived category may be included into the following short exact
sequence of groups:
$$
0\to\ZZ\oplus(A\myt\hat A)_k\to\Auteq\db A\to U(A\myt\hat A)\to1.
$$
\endproclaim

Let us study the group $\Auteq\db A$ in more detail. It has a normal
subgroup $(A\myt\hat A)_k$ consisting of the functors
$T_{a*}(\cdot)\ot P_\alpha$, where $(a,\alpha)\in A\myt\hat A$.
The quotient with respect to this subgroup is a central extension of
the group $U(A\myt\hat A)$ by~$\ZZ$. We denote this central extension by
$\myw U(A\myt\hat A)$. There are short exact sequences
$$
\gather
0\to(A\myt\hat A)_k\to\Auteq\db A\to\myw U(A\myt\hat A)\to1,
\tag4.4
\\
0\to\ZZ\to\myw U(A\myt\hat A)\to U(A\myt\hat A)\to1.
\tag4.5
\endgather
$$
To describe the central extension~\thetag{4.5}, it suffices to present
the corresponding 2-cocycle $\lambda(g_1,g_2)\in\ZZ$, where
$g_1$,~$g_2$ are elements of $U(A\myt\hat A)$.

By Proposition 3.2, if $\Phi_{\sE}$ is an equivalence, then the object
$\sE\in\db{A\myt A}$ is isomorphic to $E[k]$ for some sheaf~$E$ on
$A\myt A$. Let $E$,~$F$ be sheaves on $A\myt A$ which determine
autoequivalences~$\Phi_E$,~$\Phi_F$. Then the composite
$\Phi_F\circ\Phi_E$ is presented by some object~$G[k]$, where
$G$ is a sheaf. We put $\lambda(f_F,f_E)=k$. It is clear that $\lambda$
is a 2-cocycle, and it determines the central extension~\thetag{4.5}.

Let us calculate $\lambda$ in terms of elements of $U(A\myt\hat A)$. We put
$k=\CC$ for simplicity. Then $A\cong\CC^n/\ZZ^{2n}$, and every line bundle
$L$ on~$A$ determines a Hermitian form~$H(L)$ on~$\CC^n$
(see~\cite{11}). We denote by~$p(H)$ the number of positive eigenvalues
of~$H$. This yields a function $p\:\NS(A)\to\ZZ$. It may be extended to
$\NS(A)\ot\RR$ by
$$
p\Bigl(\,\sum_ir_i[L_i]\Bigr)=p\Bigl(\,\sum_ir_iH(L_i)\Bigr).
$$
We thus get a lower semicontinuous function $p$ on the whole of
$\NS(A)\ot\RR$. (It is easy to define the function~$p$ for any algebraically
closed field. Indeed, take~$p(L)$ to be equal to the number of negative roots
of the polynomial $P(n)=\chi(L\ot M^n)$, where $M$ is some ample line bundle.
Then extend $p$ to the whole of $\NS(A)\ot\RR$ in the way described above.)

We now consider any two elements of $U(A\myt\hat A)$,
$$
g_1=\pmatrix
x_1 & y_1
\\
z_1 & w_1
\endpmatrix, \qquad
g_2=\pmatrix
x_2 & y_2
\\
z_2 & w_2
\endpmatrix
\tag4.6
$$
such that $y_1$ and $y_2$ are isogenies. This means that there are inverse
elements~$y_1^{-1}, y_2^{-1}\in\Hom(A,\hat A)\ot\QQ$. We introduce the
notation
$$
g_1g_2=\pmatrix
x_3 & y_3
\\
z_3 & w_3
\endpmatrix.
$$
Consider the element $y_1^{-1}y_3y_2^{-1}=y_1^{-1}x_1+w_2y_2^{-1}$ of the
group $\Hom(A,\hat A)\ot\QQ$. We easily see that it belongs to
$\NS(A)\ot\QQ$. It is easy to prove the equation
$$
\lambda(g_1,g_2)=p(y_1^{-1}y_3y_2^{-1})-n.
$$
This yields a formula for $\lambda(g_1,g_2)$ in the case when
$y_1$ and~$y_2$ are invertible. Since~$\lambda$ is a cocycle, it is
uniquely determined by this formula and is uniquely extended to the
whole group $U(A\myt\hat A,\RR)\subset\EEnd(A\myt\hat A)\ot\RR$.

\example{Example 4.15}
Consider the example $A=E^n$, where $E$ is an elliptic curve without
complex multiplication. Then the group $U(A\myt\hat A)$ is isomorphic to
$\Sp_{2n}(\ZZ)$. It is well known that the fundamental group of the real
symplectic group $G=\Sp_{2n}(\RR)$ is isomorphic to~$\ZZ$ and there is a
universal central extension~$\myw G$. Moreover, the symplectic group
carries a $\ZZ$-valued 2-cocycle
$$
\mu(g_1,g_2)=\tau(l,g_1l,g_1g_2l),
$$
where $\tau$ is the Maslov index and $l$ is a Lagrangian subspace
(see, for example,~\cite{7}). There is a formula for the cocycle~$\mu$
in the matrix notation similar to~\thetag{4.6}. In the previous notation
it takes the form
$$
\mu(g_1,g_2)=\sign(y_1^{-1}y_3y_2^{-1}),
$$
the right-hand side being the signature of the symmetric matrix
$y_1^{-1}y_3y_2^{-1}$.

Comparing the cocycles $\lambda$ and $\mu$, we easily see that the cocycle
$(2\lambda-\mu)$ is trivial as an element of the second cohomology group.
Moreover, it is known that the second cohomology group of the symplectic
group $G=\Sp_{2n}(\RR)$ is equal to~$\ZZ$, and the generator determines the
universal extension~$\myw G$. It is also known that the Maslov cocycle
$\mu$ equals four times the generator. Hence $\myw G_\lambda$ is included
into an exact sequence $1\to\myw G\to\myw G_\lambda\to\ZZ/2\ZZ\to0$,
which actually splits, that is, $\myw G_\lambda$ is isomorphic to
$\myw G\myt\ZZ/2\ZZ$.
\endexample

\example{Example 4.16}
Consider an abelian variety~$A$ with the endomorphism ring $\EEnd(A)=\ZZ$.
Then the N\'eron--Severi group~$\NS(A)$ is isomorphic to~$\ZZ$. We denote
by~$L$ and~$M$ the generators of $\NS(A)$  and $\NS(\hat A)$ respectively.
The composite $\phi_M\circ\phi_L$ is equal to $N\id_A$ for some $N>0$.
The group $U(A\myt\hat A)$ coincides with the congruence subgroup
$\Gamma_0(N)\subset\SL(2,\ZZ)$. Furthermore, let $B$ be another abelian
variety such that $B\myt\wht B$ is isomorphic to $A\myt\hat A$. It is easily
verified that any such isomorphism is isometric. The abelian variety $B$
can be presented as the image of some morphism
$A@>(k\cdot\id,m\phi_L)>>A\myt\hat A$, and we may assume that
$\roman{GCD}(k,m)=1$. We denote this morphism from~$A$ to~$B$ by~$\psi$.
The kernel of~$\psi$ equals $\Ker(m\phi_L)\cap A_k$. Since
$\roman{GCD}(k,m)=1$, we actually have $\Ker(\psi)=\Ker(\phi_L)\cap A_k$.
On the other hand, we know that $\Ker(\phi)\subset A_N$. Thus we may assume
without loss of generality that $k$ divides~$N$. Each divisor~$k$ of~$N$
induces an abelian variety $B:=A/\bigl(\Ker(\phi_L)\cap A_k\bigr)$, and it
is clear that different divisors of~$N$ induce non-isomorphic abelian
varieties. Moreover, it is easy to verify that the embedding of~$B$ into
$A\myt\hat A$ splits if and only if $\roman{GCD}(k,N/k)=1$. Hence the
number of abelian varieties~$B$ with $\db B\simeq\db A$
is equal to~$2^s$,  where $s$ is the number of prime divisors of~$N$.

We now additionally assume that the abelian variety $A$ is principally
polarized, that is, we have $N=1$. In this case, if
$\db A\simeq\db B$, then $B\cong A$. Moreover, the group
$U(A\myt\hat A)$ is isomorphic to $\SL(2,\ZZ)$. As proved in~\cite{14},
the sequence~\thetag{4.4} splits for a principally polarized abelian variety.
Hence we get a description of $\Auteq\db A$ as a semidirect product of the
normal subgroup $(A\myt\hat A)_k$ and the group $\myw U(A\myt\hat A)$:
$$
\Auteq\db A\cong\myw U(A\myt\hat A)\ltimes(A\myt\hat A)_k,
$$
and $\myw U(A\myt\hat A)$ is the central extension of $\SL(2,\ZZ)$
given by a cocycle~$\lambda$. One can show that there is a sequence
$$
1\to B_3\to\myw U(A\myt\hat A)\to\ZZ/2n\ZZ\to0,
$$
where $B_3$ is the braid group of 3 strings. (We recall that $B_3$ is the
universal central extension of~$\SL(2,\ZZ)$, which is also induced by the
universal covering~$\SL(2,\RR)$.)
\endexample

\head
Appendix
\endhead

This Appendix is devoted to the proof of the statement which was made after Proposition 4.12 (see also Lemma 4.1 from \cite{14}).
\proclaim{Statement}
Any isometric isomorphism $f\:A\myt\wht A\to B\myt\wht B$ can be represented as a composition of two maps
$f_1\in U(A\myt\wht A,B\myt\wht B)$ and $f_2\in U(A\myt\wht A)$
such that $y_1\: \wht A\to B$ and $y_2\: \wht A\to A$ are isogenies.
\endproclaim
\demo{Proof}
Let us consider an isometric isomorphism $f\:A\myt\wht A\to B\myt\wht B.$ As in $\S\, 4$, we represent $f,\; f^{-1},$ and $\id_{B\times \wht B}= f \cdot f^{-1}$ in the matrix form:
$$
f=\pmatrix
x & y
\\
z & w
\endpmatrix,
\quad
f^{-1}=
\pmatrix
\wht w & - \wht y
\\
- \wht z & \wht x
\endpmatrix,
\quad
f\cdot f^{-1}=
\pmatrix
x\wht w - y\wht z & x\wht y - y\wht x
\\
z\wht w - w\wht z& w\wht x- y\wht z
\endpmatrix.
$$

Suppose that $y\:\wht A\to B$ is not an isogeny.
Let $K=\Ker y$ be an abelian subvariety that is the kernel of $y.$ Denote by $i\: K\hookrightarrow \wht A$ the embedding of the kernel of $y$ to $\wht A.$
Let us take an abelian subvariety $j\: Z\hookrightarrow\wht A$ that is a complement to $K$ in $\wht A,$ i.e. the intersection $Z\cap K$ is finite and $Z+K=\wht A.$ Consider the projection
$p\:\wht A\to K'=\wht A /Z.$ The composition $p\, i\: K\to K'$ is an isogeny. Let $M'$ be an ample line bundle on $K'.$ It induces a line bundle
$M=p^*M'$ on $\wht A$ such that the restriction of $M$ to $K$ is  also ample while the restriction to $Z$ is trivial. The line bundle $M$ gives a map
$\phi_{M}:\wht A\to A$ such that $\Ker\phi_{M}= Z$ and $\Im \phi_{M}=\wht K'\subset A.$ This implies that $\Im\phi_{M}=\Im\phi_M i$ and the composition
$\hat i\, \phi_{M} i\: K\to \wht K$ is an isogeny. The map $\phi_{M}$ defines an isometric isomorphism
$$
g_M=\pmatrix
1 & \phi_M
\\
0 & 1
\endpmatrix.
$$
from $A\myt\wht A$ to itself. Take the composition $f'= f\cdot g_M:$
$$
f'=f\cdot g_M
=\pmatrix
x & y'
\\
z & w'
\endpmatrix
=\pmatrix
x & x\phi_M+ y
\\
z & z\phi_M + w
\endpmatrix.
$$

Let us show that $y'=x\phi_M+ y$ is an isogeny. It is equivalent to check that $\Im y'\subseteq B$ coincides with the whole variety $B.$
Since the intersection $Z\cap K$ is finite and $Z+K=\wht A,$ the restriction of $y$ to $Z$ induces an isogeny $y\, j\: Z\to \Im y\subset B,$ i.e. $\Im y=\Im y\, j.$
Since $\Ker\phi_{M}= Z,$ we get $y'j=y\, j.$ Thus, we have $\Im y'\supseteq\Im y.$

Consider the maps $\wht y\: \wht B\to A$ and $\phi_M\:\wht A\to A.$ We know that $\Im\phi_M$ is equal to $\wht K'\subset A$ while $\Im\wht y$
coincides with the kernel $\Ker \wht i$ of the projection $\wht i\: A\to \wht K.$
Hence, the intersection $\Im\phi_M\cap \Im \wht y$ is finite and $\Im\phi_{M}+\Im \wht y=A.$
We also know that $\Im x+\Im y=B.$ Now, we have to use the property $x \wht y= y \wht x.$ It implies
$$
B=\Im x+ \Im y=\Im x\phi_M + \Im x\wht y+\Im y=\Im x\phi_M + \Im y \wht x+\Im y=\Im x\phi_M +\Im y.
$$
On the other hand, we know that $\Im x\phi_M=\Im x\phi_M i$ and $\Im y=\Im y\, j.$
Moreover, we have $x\phi_M i=y'i$ and $y\, j=y'j.$ Thus, we get
$$
B=\Im x\phi_M +\Im y=\Im x\phi_M i +\Im y\, j=\Im y'i +\Im y'j=\Im y'.
\nopagebreak
$$
This means that the map $y'\:\wht A\to B$ is an isogeny.

Now, any line bundle $M$ can be represented as a tensor product $L_1\otimes L_2^{-1},$ where $L_1, L_2$ are very ample line bundles on $\wht A.$
This implies the following equality for maps: $\phi_M=\phi_{L_1}-\phi_{L_2},$ where $\phi_{L_1}, \phi_{L_2}$ are isogenies.
Thus, the isometric isomorphism $f=f'\cdot g_{M}^{-1}$ can be presented as a product $f'\cdot g_{L_2}\cdot g_{L_1}^{-1} $ of three maps 
$f'\in U(A\myt\wht A,B\myt\wht B)$ and $g_{L_2}, g_{L_1}^{-1}\in U(A\myt\wht A)$
such that their $y$--components
$y',  \phi_{L_2}, -\phi_{L_1} $ are isogenies.

Finally, if we want to obtain $f$ as a composition of only two elements, then, replacing $L_2$  if necessary by its greater power $L_2^{\otimes k},$ we can ensure that the $y$-component of the product  $f_1=f'\cdot g_{L_2}$ is an isogeny as well.
Thus, we can take $f_1=f'\cdot g_{L_2}$ and $f_2=g_{L_1}^{-1}.$
\enddemo


\enddocument